\documentclass{article}
\usepackage[utf8]{inputenc}

\usepackage{authblk} 

\title{\textbf{Relax-and-round strategies for solving the Unit Commitment problem with AC Power Flow constraints}}

\author[1,2]{Dolores Gómez}
\author[3]{Simone Göttlich}
\author[1,2]{Alfredo Ríos-Alborés}
\author[1,2]{Pilar Salgado}

\affil[1]{Universidade de Santiago de Compostela, Departamento de Matemática Aplicada, Santiago de Compostela, 15782, Spain}
\affil[2]{Galician Centre for Mathematical Research and Technology (CITMAga), Santiago de Compostela, 15782, Spain}
\affil[3]{University of Mannheim, School of Business Informatics and Mathematics, 68159 Mannheim, Germany}

\date{}

\usepackage[numbers]{natbib}
\usepackage{amssymb}
\usepackage{amsfonts,amsmath,pdfsync,color,graphicx}
\usepackage{amsthm}
\usepackage{ upgreek }
\usepackage{afterpage}
\usepackage{graphicx}
\usepackage{geometry}
\usepackage{enumitem}
\usepackage{wrapfig}
\usepackage{subcaption} 
\usepackage{algorithm}
\usepackage{algpseudocode}
\usepackage{longtable}
\usepackage{booktabs}

\usepackage{pifont}
\usepackage{xcolor}

\newcommand{\cmark}{\textcolor{green}{\ding{51}}} 
\newcommand{\xmark}{\textcolor{red}{\ding{55}}}   

\definecolor{dkgreen}{rgb}{0,0.6,0}
\definecolor{gray}{rgb}{0.5,0.5,0.5}
\definecolor{mauve}{rgb}{0.58,0,0.82}

\usepackage{bm}

\usepackage{listings}
\usepackage{relsize}
\usepackage{multirow} 

\makeatletter
\renewcommand*\env@matrix[1][\arraystretch]{%
	\edef\arraystretch{#1}%
	\hskip -\arraycolsep
	\let\@ifnextchar\new@ifnextchar
	\array{*\c@MaxMatrixCols c}}
\makeatother

\newcommand{\tv}{\tilde{v}}
\newcommand{\tu}{\tilde{u}}
\newcommand{\tw}{\tilde{w}}
\newcommand{\tP}{\tilde{P}}
\newcommand{\tQ}{\tilde{Q}}
\newcommand{\defeq}{\mathrel{\mathop:}=}

\begin{document}
\maketitle

\begin{center}
\begin{longtable*}[tbp]{rl}
\multicolumn{2}{c}%
{\small{\textbf{Nomenclature}}} \\[1ex]
\multicolumn{1}{c}{\textrm{Symbol}} & \multicolumn{1}{c}{\textrm{Name and description}}  \\ \midrule
\endfirsthead

\multicolumn{2}{c}%
{\small{\textbf{Nomenclature}}} \\[1ex]
 \multicolumn{1}{c}{\textrm{Symbol}} & \multicolumn{1}{c}{\textrm{Name and description}} \\ \midrule
\endhead

 \multicolumn{2}{r}{\small{Continue on next page}} \\ \midrule
\endfoot

\bottomrule
\endlastfoot

 $a_{g,2},a_{g,1},a_{g,0}$     & coefficients of the quadratic cost function for generation unit $g$\\[0.7ex]
             $ a_{il} $      &  coefficient $il$ of matrix $\mathcal{A}$\\[0.5ex]
             $\mathcal{A}$      & incidence matrix of $\mathcal{P}$\\[0.5ex]
             $b_{ik}$      & imaginary part of coefficient $ik$ of admittance matrix $\mathbb{Y}$ \\[0.5ex]
             $c_{g,u},c_{g,d}$      & start-up/shut-down constant costs of generation unit $g$\\[0.5ex]
             $ c_{1,il},c_{2,il}$      & coefficient $il$ of matrices $C_1$,$C_2$\\[0.5ex]
             $ C_1,C_2$      & auxiliary matrices defined from $\mathcal{A}$\\[0.5ex]
             $g $      & generation unit index\\[0.5ex]
             $g_{ik}$      & real part of coefficient $ik$ of admittance matrix $\mathbb{Y}$\\[0.5ex]
             $G $      &  set of generation units\\
             $ G_i$      & set of generation units placed at bus $i$\\[0.5ex]
             $ G_r$      & set of generation units in region $r$\\
             $ G_{\text{th}}$      & set of thermal generation units\\[0.5ex]
             $ G_{\text{sc}}$      & set of synchronous condensers\\[0.5ex]
             $ i$      &  bus index\\
             $ I $      & set of buses\\[0.5ex]
             $ I_D$      & set of buses where power is demanded\\
             $ l $      & line index\\[0.5ex]
             $ L $      & set of lines\\[0.5ex]
             $ N_{\text{b}} $ & number of buses \\[0.5ex]
             $ N_{\text{l}} $ & number of elements \\[0.5ex]
             $ \mathcal{N}\left(i\right)$      & set of bus $i$ neighbours\\[0.5ex]
             $\mathcal{P} $ &  directed graph modeling the power grid\\[0.5ex]
             $ P_{g,t}^G$      & active power generated by power plant $g$ at time $t$\\[0.5ex]
             $ \tP_{g,t}^G$      & relaxation of  $ P_{g,t}^G$ to $\left[0,P_{\max,g}\right]$\\[0.5ex]
             $ P_{i,t}^D$      & active power demand at bus $i$ at time $t$\\[0.5ex]
             $ P_{\max,g}, P_{\min,g}$      & maximum/minimum active power generation limits for power plant $g$\\[0.5ex]
             $P_{g,t}^{\text{res}}$ & active power reserve for generating unit $g$ at time $t$\\[0.5ex]
             $P_{r,t}^{\text{res}}$ & active power reserve requirements for region $r$ at time $t$\\[0.5ex]
             $ Q_{g,t}^G$      & reactive power generated by power plant $g$ at time $t$\\[0.5ex]
             $ \tQ_{g,t}^G$      &  relaxation of  $ Q_{g,t}^G$ to $\left[\min\{Q_{\min,g}, 0\}, \max\{Q_{\max,g},0\}\right]$\\[0.5ex]
             $ Q_{i,t}^D$      & reactive power demand at bus $i$ at time $t$\\[0.5ex]
             $ Q_{\max,g}, Q_{\min,g}$      & maximum/minimum reactive power generation limits for power plant $g$\\[0.5ex]
             $ r$      &  region index\\
             $R $      &  set of regions\\[0.5ex]
             $R_{g}^{d}, R_{g}^{u}$      & active power ramp down/up  specifications for generation unit $g$\\
             $S_{g}^{d}, S_{g}^{u}$      & shutting-down/starting-up active power generation limits of power plant $g$\\
             $\vec{\mathbb{S}}$      & vector of complex power injections from the buses to the grid\\[0.5ex]
             $t $      & time index\\[0.5ex]
             $ T$      &  set of discrete times\\
             $ T_{g}^{d},T_{g}^{u}$      & minimum down/up times for generation unit $g$\\[0.5ex]
              $u_{g,t} $      & on/off status of generation unit $g$ at time $t$\\[0.5ex]
             $ \tu_{g,t} $      & continuous relaxation of $u_{g,t}$ to $\left[0,1\right]$ \\[0.5ex]
             $v_{g,t} $      & turn on decision variable for generation unit $g$ at time $t$\\[0.5ex]
             $ \tv_{g,t} $      & continuous relaxation of $v_{g,t}$ to  $\left[0,1\right]$ \\[0.5ex]
             $V_{i,t} $      & complex voltage magnitude of bus $i$ at time $t$\\[0.5ex]
             $V_{\max,i},V_{\min,i} $      & maximum/minimum complex voltage magnitude limits for bus $i$\\[0.5ex]
             $\vec{\mathbb{V}}$      & vector of complex voltage for the buses\\[0.5ex]
             $w_{g,t} $      & turn off decision variable for generation unit $g$ at time $t$\\[0.5ex]
             $ \tw_{g,t} $      & continuous relaxation of $w_{g,t} $ to $\left[0,1\right]$ \\[0.5ex]
             $\mathbb{Y}$      & admittance matrix of the power grid\\[0.5ex]
             $\mathbb{Y}_1,\mathbb{Y}_2$      & from/to components of the admittance matrix $\mathbb{Y}$\\[0.5ex]
              $\vec{\mathbb{Y}}_{\text{sh}}$      & complex vector of bus shunt elements\\[0.5ex]
             $\theta_{i,t}$      & complex voltage angle  of bus $b$ at time $t$\\
\end{longtable*}
\end{center}
\section{Abstract}
The Unit Commitment problem with AC power flow constraints (UC-ACOPF) is a non-convex mixed-integer nonlinear programming (MINLP) problem encountered in power systems. Its combination of combinatorial complexity and non-convex nonlinear constraints makes it particularly challenging. A common approach to tackle this issue is to relax the integrality condition, but this often results in infeasible solutions. Consequently, rounding heuristics are frequently employed to restore integer feasibility.

This paper addresses recent advancements in heuristics aimed at quickly obtaining feasible solutions for the UC-ACOPF problem, focusing specifically on direct relax-and-round strategies. We propose a model-based heuristic that rescales the solution of the integer-relaxed problem before rounding. Furthermore, we introduce rounding formulas designed to enforce combinatorial constraints and aim to maintain AC feasibility in the resulting solutions. These methodologies are compared against standard direct rounding techniques in the literature, applied to a 6-bus and a 118-bus test systems. Additionally, we integrate the proposed heuristics into an implementation of the Feasibility Pump (FP) method, demonstrating their utility and potential to enhance existing rounding strategies.

\section{Introduction} \label{Sect0}

The Unit Commitment (UC) problem is a fundamental optimization task in managing and planning the operation of power grids. It aims to determine an optimal schedule for power plant operations ensuring generation meets demand \cite{garver1962power} complying with the technical and operational constraints of the different power generation technologies, as well as with the power system and the transmission characteristics of the grid \cite{padhy2004unit}. Over the years, the UC problem evolved in response to different power market paradigms, transitioning from centralized to competitive markets, and adapting to changes in the energy mix, which now, for example, includes a growing share of renewable energy sources alongside traditional coal, gas, hydro, and nuclear power plants. Additionally, the UC problem can be used as a short-term real-time tool \cite{li2024objective}, \cite{wang2013stochastic}, or as a long-term planner \cite{zhang2023mathematical}, but generally it is used as a day-ahead scheduler.

From a mathematical point of view, the UC problem can be formulated in different ways depending on the purpose. For instance, it can be stated as a cost minimization or as a profit maximization problem. Generation costs are mainly driven by fuel consumption that is often modeled by quadratic functions or their piecewise linear approximations \cite{frangioni2008tighter}, but other factors such as water reserves \cite{caroe1998two} and uncertain weather conditions \cite{wang2008security}, \cite{wang2011wind}, can also play a significant role. The operational characteristics of power plants over time are represented by linear constraints that limit ramp-up/ramp-down rates or enforce minimum on/off times. The UC problem is typically a mixed-integer problem because it involves binary variables that model the on/off status of power plants. Available formulations include the well-known 1-bin \cite{carrion2006computationally} and 3-bin \cite{morales2013tight}, and alternative ones \cite{yang2021two}, \cite{atakan2017state}. Though, other modeling approaches avoid integer variables \cite{chen2019modeling}. Transmission capabilities are sometimes neglected, by just enforcing a global linear power balance \cite{caroe1998two}, or modeled using linear approximations such as the DC power flow equations \cite{wu2010accelerating}, or, more accurately, with the non-linear and non-convex AC power flow equations \cite{castillo2016unit}. Additionally, the spinning reserve requirements for the system should be taken into account. To summarize, depending on the UC formulation selected, the resulting problem may be a deterministic or stochastic mixed-integer linear programming (MILP) problem, mixed-integer quadratic programming (MIQP) problem, non-convex mixed-integer non-linear programming (MINLP) problem, among others.

The large size of real power grids, the combinatorial nature of the UC problem, the uncertainty in some input data, and the non-convexity and non-linearity of certain grid representations present significant challenges. In practice, however, these factors are rarely addressed all at once. Instead, the UC problem is often simplified based on the specific goals of the study; although the ultimate goal is to develop the most precise and robust optimization tools, as seen in the recent Grid Optimization Competition \cite{holzer2024grid}. Various approaches have been used to tackle the UC problem, including reformulations \cite{chen2019modeling}, \cite{yang2014tight}, \cite{ostrowski2011tight}, relaxation and/or decomposition techniques (e.g., semi-definite programming \cite{bai2009semi}, Lagrange Relaxation \cite{murillo2000parallel}, Benders' Decomposition \cite{niknam2009new}, \cite{nick2015security}, Outer Approximation \cite{castillo2016unit}), heuristic methods \cite{cobos2018network}, \cite{senjyu2003fast}, or machine learning \cite{xavier2021learning}. Improving the efficiency of MILP formulations for this problem remains an open area of research \cite{tejada2019unit2}. However, relying solely on linear grid representations for power grids can lead to additional costs due to uplift payments \cite{sauer2014uplift} or even to security issues in highly constrained grids \cite{bai2015decomposition}. This has motivated efforts to develop convex relaxations approaches for the AC power flow equations, providing more tractable formulations such as convex second-order conic relaxations \cite{bai2015decomposition}, \cite{liu2018global}, \cite{tuncer2022misocp}. Works \cite{montero2022review} and \cite{aharwar2023unit} are recent reviews on the general UC problem and different stochastic modeling frameworks, respectively.

As the combinatorial nature of the UC problem is the main challenge when solving it, a common approach is to relax the integrality condition on binary variables. State-of-the-art methods often rely on branch-and-cut solvers, where integer relaxation is integrated into the algorithm. Moreover, linear programming (LP) and non-linear programming (NLP) problems are generally more efficient to solve than their mixed-integer counterparts, making it tempting to relax mixed-integer constraints into continuous ones. For MILP formulations, work~\cite{wuijts2024effect} empirically concludes that, in certain contexts, relaxed approximations can yield significant computational time savings with minimal loss of solution quality. However, when the integrality condition is relaxed, the resulting solutions are often infeasible for the original problem. To restore integer feasibility, rigorous iterative methods, such as branching, or rounding heuristics, are typically employed.

In the literature, rounding is a simple yet useful technique often applied directly, for example, by using a basic threshold (e.g., 0.5), or through more advanced strategies involving optimization and iterative algorithms. These heuristics do not guarantee general feasibility or optimality for the original problem. In the latter category, the Feasibility Pump \cite{fischetti2005feasibility} and local integer search methods, such as Variable Neighborhood \cite{lazic2016variable}, along with their variants, are popular choices for relax-and-round strategies in mixed-integer programming (MIP), including UC problems \cite{li2024objective}, \cite{ma2020unit}. However, direct rounding techniques are also employed as straightforward strategies to ensure integer feasibility. In stochastic problems, different scenarios are often considered, and the weighted average commitment is simply rounded, as in \cite{caroe1998two}. This approach is also common in the context of machine learning, where the average commitment is interpreted as a confidence level (probability of commitment) after the training stage \cite{gao4718365topology}. Direct rounding is frequently used as a heuristic to enhance the convergence rate of progressive hedging algorithms \cite{ordoudis2015stochastic}. Procedures involving relaxation approaches usually rely on rounding at some stage to ensure the integrality of the solution \cite{wu2023synergistic}, \cite{kjeldsen2012heuristic}. Moreover, direct rounding heuristics can be more sophisticated, aiming to avoid combinatorial infeasibility related to other constraints, such as the minimum up/down time constraints \cite{yang2014tight}, \cite{bai2009semi}.

In this paper, we propose new different relax-and-round approaches involving rounding formulas for the deterministic Unit Commitment problem with AC power flow constraints (UC-ACOPF) involving mainly thermal power plants. First, we introduce a rescaling technique for the integrality-relaxed solution; it consists of a model-based heuristic aiming to improve the effectiveness of other existing rounding methods. Then, a rounding formula for a relax-and-round direct approach is proposed. In addition, all the proposed heuristics are compared with other simpler rounding techniques from the literature and validated using two common benchmarks: the 6-bus test system and a modified IEEE 118-bus system. Finally, the effectiveness of the rescaled technique on other rounding heuristics is tested by implementing a Feasibility Pump algorithm. The  paper is organized as follows: Section \ref{Sect1} introduces the Unit Commitment problem. In Section \ref{Sect2}, we present direct relax-and-round heuristics, followed by a discussion of their numerical results in Section \ref{Sect3}. Section \ref{Sect5} covers the implementation of a Feasibility Pump algorithm combined with our proposed relax-and-round heuristics, along with the corresponding numerical results.

\section{The UC-ACOPF problem}\label{Sect1}

Let $\mathcal{P}= \left(I, L\right)$ be a directed graph modeling the topology of a power grid, where nodes $I=\{1, \dots, N_{\text{b}}\}$ represent the buses, and edges $L\subset I\times I$ the different $N_{\text{l}}$ elements, like transmission lines, transformers, etc. Then, for all $l\in L$, $l=\left(i,k\right)$ for some $i,k\in I$, and we say that $i$ and $k$ are neighbors or connected in $\mathcal{P}$. Let $\mathcal{N}(i)$ represent the set of indices corresponding to the grid neighbors of node $i$.

As we are dealing with a directed graph, some of the power grid linear relationships between the problem variables can be deduced, defined, or expressed using matrices. For this purpose, it is useful to define the incidence matrix of the graph, namely, $\mathcal{A}:=\left(a_{kl}\right)\in\mathcal{M}_{N_{\text{b}}\times N_{\text{l}}}$, where
$$ a_{kl} =
\left\{
\begin{array}{rl}
     0, & \text{if node } k \text{ not belong to line } l, \\
    -1, & \text{if node } k \text{ is the `\text{from}' node of line } l, \\
     1, & \text{if node } k \text{ is the `\text{to}' node of line } l.
\end{array}
\right.$$
Also, we define some auxiliary matrices: $C_1 := \left(c_{1,lk}\right)\in\mathcal{M}_{N_{\text{l}}\times N_{\text{b}}}$, where $c_{1,lk} = \max\{0,-a_{kl}\}$, and $C_2 := \left(c_{2,lk}\right)\in\mathcal{M}_{N_{\text{l}}\times N_{\text{b}}}$, where $c_{2,lk} = \max\{0,a_{kl}\}$.

Placed at some buses there are generation units $g\in G=\{1,\dots, N_{\text{g}}\}$. Let $G_i\subset G$ denote the subset of generation units placed at node $i\in I$ in the graph. Also, the nodes are grouped into regions $r\in R$, a subset of buses and elements of the power grid grouped for technical reasons; let $G_r\subset G$ denote the subset of generation units placed at nodes in region $r$.

Let us consider an initial time $t_0$, and the set of discrete times $T = \{t_1,\dots,t_f\}$. The decision values of the variables at $t_0$ are known. There is also a known demand for active $P_{i,t}^D$ and reactive $Q_{i,t}^D$ power at buses $i\in I_D\subset I$, $t\in T$.

Classically, the Unit Commitment (UC) problem is formulated as a mixed-integer linear constrained problem, where the cost function often includes quadratic terms, sometimes approximated by piecewise linear functions. The mixed-integer linear constraints capture the on/off logic of power plants, ramp-up/ramp-down limits, minimum up/down time requirements, and power generation limits. This formulation is sometimes referred to as the `Unit Commitment skeleton' (e.g., \cite{liu2018global}) a term we also adopt in this work. Typically, a global linear power balance equation is incorporated into the UC skeleton. However, in this work, we model the power grid using the AC power flow (ACPF) equations. This introduces non-linear constraints, including the active and reactive power balances at each bus, and the complex current limits of the transmission lines, which are enforced via the ACPF formulation. We will distinguish between different levels of feasibility depending on whether a point, comprising the values of all decision variables, satisfies only the UC skeleton constraints, only the ACPF constraints, or both.

In what follows, we will refer to edges as lines. Although not every element in the power grid is strictly a transmission line, all components will be modeled using the same linear circuit for the ACPF equations, the so-called $\pi$-model \cite{zimmerman2024matpower}. Similarly, we will use the terms nodes and buses indistinctly. Also, we now introduce a general notation for the decision variables: let $\boldsymbol{x}$ denote the vector of all continuous decision variables for all time steps $t$ or other sub-indices (e.g., $g\in G$, $i\in I$, etc.), and let $\boldsymbol{y}$ denote the vector of binary decision variables. Thus, a pair $\left(\boldsymbol{x},\boldsymbol{y}\right)$ represents a point in the space of decision variables for the UC-ACOPF problem.

\subsection{Unit Commitment skeleton}

In this paper, for simplicity, we consider only thermal power plants, $g\in G_{\text{th}}\subset G$, and also synchronous condenser units for reactive power compensation, $g\in G_{\text{sc}}\subset G$. Nuclear power plants or renewable units, such as windmill fields or hydro-power plants, could also be considered. In the literature, there are several equivalent formulations available for modeling the UC skeleton constraints, and each has potential advantages in terms of computational efficiency \cite{carrion2006computationally}, \cite{ostrowski2011tight}. This has been a hot topic in the literature for some time; a well-known formulation including a discussion about compactness and tightness as efficiency factors for MILP problems was presented in \cite{morales2013tight}. We choose to combine the formulations from \cite{liu2018global} and \cite{constante2022ac}, as in this study we focus on the performance of our approach in terms of the feasibility properties of the returned solutions. However, since the proposed techniques rely on the relaxed problem quality, the overall performance could be improved by choosing other formulations.

For each thermal generation unit $g\in G_{\text{th}}$ and time $t\in T$, we need to consider the following decision variables. The on/off status and turn on/off decisions, that are modeled using the binary variables $u_{g,t}$, $v_{g,t}$, $w_{g,t}\in \{0,1\}$, respectively; and the active $P_{g,t}^G\in \{0\}\cup\left[P_{\min,g},P_{\max,g}\right]$ and reactive $Q_{g,t}^G\in \{0\}\cup \left[Q_{\min,g},Q_{\max,g}\right]$ power generated.

The total cost function of the power system operation depends on the active power generation and the on/off scheduling decisions and is given by the quantity
\begin{equation}
\label{UC_cost}
 f_{\text{UC}}\left(\boldsymbol{x},\boldsymbol{y}\right)=\sum_{t\in T}\sum_{g\in G_{\text{th}}} a_{g,2}\left(P_{g,t}^G\right)^2+a_{g,1}P_{g,t}^G+a_{g,0}u_{g,t}+c_{g,u}v_{g,t}+c_{g,d}w_{g,t}.
\end{equation}
The UC aims to minimize (\ref{UC_cost}) subjected to the operation linear constraints presented next.

For each $g\in G_{\text{th}}$, we consider the combinatorial constraints modeling the on/off logic
\begin{equation}
\label{UC_BVL}
 u_{g,t-1} - u_{g,t} + v_{g,t} - w_{g,t} = 0,\quad t\in T,
\end{equation}
the minimum up times, enforced by
\begin{equation}
\label{UC_MTup}
 \sum_{\tilde{t}=t-T_g^u+1}^{t} v_{g,\tilde{t}} \le u_{g,t},\quad t\in\{T_g^u,\dots,t_f\},
\end{equation}
and the minimum down times, by
\begin{equation}
\label{UC_MTdown}
 \sum_{\tilde{t}=t-T_g^d+1}^{t} w_{g,\tilde{t}} \le 1 - u_{g,t},\quad t\in\{T_g^d,\dots,t_f\}.
\end{equation}
The parameters $T_g^u$ and $T_g^d$ represent the minimum time, expressed in hours, that each thermal generation has to remain up or down after a start-up or shut-down decision, respectively, for each thermal generator $g$.

Additionally, for every region $r\in R$ and time $t\in T$, there is an active power reserve requirement, $P_{r,t}^{\text{res}}$. Then, an active power reserve variable is defined for each thermal plant and time, $P_{g,t}^{\text{res}}\in\left[0,P_{\max,g}-P_{\min,g}\right]$. The active power reserve requirement constraints are expressed by the linear inequalities:
\begin{equation}
\label{UC_ReReq}
\sum_{g\in G_{r}}P_{g,t}^{\text{res}} \ge P^{\text{res}}_{r,t}.
\end{equation}
For simplicity, in this work, we will assume that only one region $r$ is considered in the power grid.

For each synchronous condenser $g\in G_{\text{sc}}$, the reactive power $Q^G_{g,t}$ is limited by
\begin{equation}
\label{UC_SClim}
 Q_{\min,g} \le Q^G_{g,t} \le Q_{\max,g} .
\end{equation}
We also need to consider a set of mixed-integer linear constraints, relating binary and continuous decision variables. For each $g\in G_{\text{th}}$ and $t\in T$, we consider the ramping-up and down constraints for the active power generation,
\begin{equation}
\label{UC_Rup}
 P^G_{g,t} +P_{g,t}^{\text{res}} - P^G_{g,t-1} \le R^u_g u_{g,t-1} + S_g^u v_{g,t},
\end{equation}
\begin{equation}
\label{UC_Rdown}
 P^G_{g,t-1} - P^G_{g,t} \le R^d_g u_{g,t} + S_g^d w_{g,t},
\end{equation}
where $R_{g}^{u}$,$R_{g}^{d}$ represent the power ramping up/down capability, and $S_{g}^{u}$,$S_{g}^{d}$ the power generation limits after starting-up and before shutting down. Also, we consider the general power generation limits
\begin{equation}
\label{UC_APGLs}
 P_{\min,g} u_{g,t} \le P^G_{g,t},\quad  P^G_{g,t} +P_{g,t}^{\text{res}} \le P_{\max,g} u_{g,t},
\end{equation}
\begin{equation}
\label{UC_RPGLs}
 Q_{\min,g} u_{g,t} \le Q^G_{g,t} \le Q_{\max,g} u_{g,t}.
\end{equation}
Usually, $P_{\min,g}>0$ defining an unconnected domain for $P_{g,t}^G$. The  power output of the generator jumps from the off state (represented by $0$) to the operating range $\left[P_{\min,g},P_{\max,g}\right]$,  depending on the binary variable $u_{g,t}$ configuration. In contrast, $Q_{\min,g}$ is usually zero or a negative value (symmetrical situation for $Q_{\max,g}$), which means $Q_{g,t}^G$ belongs to a connected set $\left[Q_{\min,g},Q_{\max,g}\right]$. This discontinuous nature of $P_{g,t}^G$ is inconvenient \cite{atakan2017state}, so this issue is sometimes mitigated redefining the variable $P_{g,t}^G$ by splitting it into two mixed-integer components, for example, $P_{g,t}^G=P_{\min,g}u_{g,t}+p^{\Delta}_{g,t}$, where $p^{\Delta}_{g,t}\in\left[0,P_{\max,g}-P_{\min,g} \right]$.
\subsection{AC power flow equations}
To model the power balance at each $i\in I$ for $t\in T$, we will use the equations
$$\sum_{g\in G_i} P^G_{g,t} - P_{i,t}^D - \mathcal{R}\text{e}\left(\left(\vec{\mathbb{S}}\left(\vec{\mathbb{V}}_t\right)\right)_i\right)= 0, $$
$$\sum_{g\in G_i} Q^G_{g,t} - Q_{i,t}^D - \mathcal{I}\text{m}\left(\left(\vec{\mathbb{S}}\left(\vec{\mathbb{V}}_t\right)\right)_i\right)= 0, $$
where $\vec{\mathbb{S}}\left(\cdot\right)\in \mathbb{C}^{N_{\text{b}}}$ is the complex vector function returning net power injections from the bus to the grid.

The AC power flow (ACPF) model states that $\vec{\mathbb{S}}$ is a function of $\vec{\mathbb{V}} \in \mathbb{C}^{N_{\text{b}}}$, the vector of buses complex voltages at each bus, under the relationship given by
$$\vec{\mathbb{S}}\left(\vec{\mathbb{V}}\right) = \text{diag}\left(\vec{\mathbb{V}}\right)\text{conj}\left(\mathbb{Y}\vec{\mathbb{V}}\right), $$
where $\mathbb{Y}\in\mathcal{M}_{N_{\text{b}}\times N_{\text{b}}}\left(\mathbb{C}\right)$ is referred to as admittance matrix of the power grid, and $\text{diag}\left(\cdot\right)$ and $\text{conj}\left(\cdot\right)$ are, respectively, the diagonal matrix and the complex conjugation operators for vectors. We denote $\left(\mathbb{Y}\right)_{ik}=g_{ik}+jb_{ik}$, for all $i,k\in I$, being $j$ the imaginary unit.

Sometimes, instead of using the grid admittance matrix directly, the following equivalence is used:
$$ \mathbb{Y} = C_1^t\mathbb{Y}_1 + C_2^t \mathbb{Y}_2 + \text{diag}\left(\vec{\mathbb{Y}}_{\text{sh}}\right), $$
where $\vec{\mathbb{Y}}_{\text{sh}}\in \mathbb{C}^{N_{\text{b}}}$ is a complex vector, comprising the coefficients of the shunt elements at the buses. Similarly, we denote $\left(\vec{\mathbb{Y}}_{\text{sh}}\right)_i=g^{\text{sh}}_{i}+jb_{i}^{\text{sh}}$. There,
$$\mathbb{Y}_1 \defeq \text{diag}\left(\vec{\mathbb{Y}}_{L,11}\right)C_1 +
\text{diag}\left(\vec{\mathbb{Y}}_{L,12}\right)C_2,$$
$$\mathbb{Y}_2 \defeq \text{diag}\left(\vec{\mathbb{Y}}_{L,21}\right)C_1 +
\text{diag}\left(\vec{\mathbb{Y}}_{L,22}\right)C_2,$$
where $\vec{\mathbb{Y}}_{L,mn}\in \mathbb{C}^{N_{\text{l}}}$ are the complex vectors containing the $\left(m,n\right)\in \{1,2\}\times\{1,2\}$ components of the element-admittance matrix defined from the $\pi$-model, for every element in the power grid.

Furthermore, the line capacity of transmitted power or current is limited, in terms of complex power or current, by one of the following inequalities
$$
\left|\mathbb{S}_{l,t}\right|^2\le S_{\max,l}^2,\quad \text{or}\quad \left|\mathbb{I}_{l,t}\right|^2\le I_{\max,l}^2, \text{ for each } l\in L,\ t\in T,
$$
where $\mathbb{S}_{l,t}$ and $\mathbb{I}_{l,t}$ are, respectively, the apparent complex power and the complex current through line $l$ at time $t$. The `from' and `to' components of those magnitudes are taken into account. We will consider the line capacity limits in terms of the current, as in \cite{castillo2016unit}.

\subsubsection{Polar coordinates for the complex voltage}
There are different equivalent formulations for the ACPF model equations depending on the chosen complex representation for the modeling variables. We adopt the polar coordinate representation for the complex voltages. Consequently, the decision variables are the voltage magnitude, $V_{i,t}\in\left[V_{\min,i},V_{\max,i}\right]$, and the voltage phase, $\theta_{i,t}\in\left[-\pi, \pi\right]$, for each node $i\in I$ and time $t\in T$.

The bus power balance is then modeled, for all $i\in I$ and $t\in T$, by the non-convex non-linear equations
\begin{equation}
\label{ACPF_PV_APB}
\sum_{g\in G_i} P^G_{g,t} - P_{i,t}^D -
V_{i,t}\sum_{k\in\mathcal{N}\left(i\right)}V_{k,t}\left(g_{ik}\cos\left(\theta_{i,t}-\theta_{k,t}\right)+b_{ik}\sin\left(\theta_{i,t}-\theta_{k,t}\right)\right) = 0,
\end{equation}
\begin{equation}
\label{ACPF_PV_RPB}
\sum_{g\in G_i} Q^G_{g,t} - Q_{i,t}^D -
V_{i,t}\sum_{k\in\mathcal{N}\left(i\right)}V_{k,t}\left(g_{ik}\sin\left(\theta_{i,t}-\theta_{k,t}\right)-b_{ik}\cos\left(\theta_{i,t}-\theta_{k,t}\right)\right) = 0.
\end{equation}
The complex voltage magnitude at each bus is constrained by the box bounds
\begin{equation}
\label{ACPF_PV_VML}
V_{\min,i}\le V_{i,t}\le V_{\max,i},
\end{equation}
representing the technical specifications for the bus voltage in the power grid. The complex angle is mathematically limited by
\begin{equation}
\label{ACPF_PV_VAL}
-\pi \le \theta_{i,t}\le \pi,
\end{equation}
and for the reference node, the voltage angle is fixed, i.e.,
\begin{equation}
\label{ACPF_PV_VALref}
\theta_{i,t}=\theta_{i,0},\quad \forall t\in T.
\end{equation}
The current limits of line $l\in L$ for time $t\in T$ using polar coordinates can be expressed in matrix form as
\begin{equation}
\label{ACPF_PV_LineLim}
\left(\text{diag}\left(\mathbb{Y}_m\vec{\mathbb{V}}_t\right)\text{conj}\left(\mathbb{Y}_m\vec{\mathbb{V}}_t\right)\right)_l\le \left(I_{\max,l}\right)^2,
\end{equation}
with $m\in\{1,2\}$ representing the `to' and `from' components and $\left(\vec{\mathbb{V}}_t\right)_{i}=V_{i,t}\mathrm{e}^{j\theta_{i,t}}$ for each $i\in I$.

\subsection{UC-ACOPF minimization problem}
The main object of study, the unit commitment problem with AC power flow constraints is a non-convex MINLP problem, defined as the following minimization problem:
\begin{equation}
\label{UC-ACOPF_def}
\left(\text{UC-ACOPF}\right) \quad
\left\{
\begin{array}{rl}
    \min & \left(\ref{UC_cost}\right) \\[0.1ex]
   \text{s.t.:} & \left(\ref{UC_BVL}\right)-\left(\ref{ACPF_PV_LineLim}\right), \\[0.1ex]
    & V_{i,t}\in\left[V_{\min,i},V_{\max,i}\right],\ \theta_{i,t}\in\left[-\pi,\pi\right],\ i\in I, t\in T,\\[0.1ex]
    & \theta_{i,t} = \theta_{i,0},\ \forall t\in T,\ \text{ for the reference node},\\[0.1ex]
    &   Q^G_{g,t}\in\left[Q_{\min,g},Q_{\max,g}\right],\ \forall \left(g,t\right)\in G_{\text{sc}}\times T,\\[0.1ex]
    & P^G_{g,t}\in\{0\}\cup\left[P_{\min,g},P_{\max,g}\right],\ \forall \left(g,t\right)\in G_{\text{th}}\times T,\\[0.1ex]
    &   Q^G_{g,t}\in\{0\}\cup\left[Q_{\min,g},Q_{\max,g}\right],\ \forall \left(g,t\right)\in G_{\text{th}}\times T,\\[0.1ex]
    &P_{g,t}^{\text{res}}\in\left[0,P_{\max,g}- P_{\min,g}\right],\ \forall \left(g,t\right)\in G_{\text{th}}\times T, \\[0.1ex]
    & u_{g,t}, v_{g,t}, w_{g,t} \in \{0,1\},\ \forall \left(g,t\right)\in G_{\text{th}}\times T. \\[0.1ex]
\end{array}
\right.
\end{equation}
Note that we have assigned labels to specific problems discussed in this paper; thus, for clarity, we will refer to (\ref{UC-ACOPF_def}) as problem (UC-ACOPF) or simply (UC-ACOPF). The same notation criteria will be used for other labeled problems.

A point $\left(\boldsymbol{x},\boldsymbol{y}\right)$ is considered UC-feasible if satisfies constraints (\ref{UC_BVL})-(\ref{UC_MTdown}) and (\ref{UC_Rup})-(\ref{UC_Rdown}); it is considered AC-feasible if it satisfies constraints (\ref{ACPF_PV_APB})-(\ref{ACPF_PV_RPB}) and (\ref{ACPF_PV_LineLim}). For simplicity in the discussion, reserve requirements and other remaining constraints are excluded from both groups. We say that a point is feasible for the UC-ACOPF problem if it satisfies all the constraints outlined in the problem definition (\ref{UC-ACOPF_def}).

\section{Relax-and-round strategies} \label{Sect2}

In this paper, we present different heuristics procedures aimed at quickly obtaining feasible points for the UC-ACOPF problem using relax-and-round strategies.{These strategies involve relaxing the integer constraints to solve a continuous optimization problem, and then rounding the solution to satisfy the integer constraints, aiming for feasible and near-optimal results}.

The Relaxed Unit Commitment (RUC) problem is defined as the linear programming (LP) problem obtained by relaxing the integrality constraint on $u_{i,t}$, $v_{i,t}$ and $w_{i,t}$ in the UC problem, allowing them to take values in $\left[0,1\right]$. Notice that $v_{g,t}$ and $w_{g,t}$ could be considered continuous variables in $\left[0,1\right]$ since the binary logic is sufficiently enforced by constraints (\ref{UC_BVL})-(\ref{UC_MTdown}) and the fact that $u_{g,t}\in\{0,1\}$. For that reason, when discussing the relaxation of the UC problem, the focus is placed on relaxing $u_{g,t}$. However, in practice, when implementing mixed-integer algorithms, it is recommended to keep $v_{g,t}$ and $w_{g,t}$ as binary variables \cite{ostrowski2011tight}, which is what we have done. Thus, the RUC-ACOPF problem is defined as the non-convex non-linear programming (NLP) problem obtained when the integrality condition on $u_{g,t}$, $v_{g,t}$ and $w_{g,t}$ is relaxed within the UC-ACOPF formulation, namely,
\begin{equation}
\left(\text{RUC-ACOPF}\right) \quad
\left\{
\begin{array}{rl}
    \min & \left(\ref{UC_cost}\right) \\[0.1ex]
   \text{s.t.:} & \left(\ref{UC_BVL}\right)-\left(\ref{ACPF_PV_LineLim}\right), \\[0.1ex]
    & V_{i,t}\in\left[V_{\min,i},V_{\max,i}\right],\ \theta_{i,t}\in\left[-\pi,\pi\right],\ i\in I, t\in T,\\[0.1ex]
    & \theta_{i,t} = \theta_{i,0},\ \forall t\in T,\ \text{ for the reference node},\\[0.1ex]
    &   Q^G_{g,t}\in\left[Q_{\min,g},Q_{\max,g}\right],\ \forall \left(g,t\right)\in G_{\text{sc}}\times T,\\[0.1ex]
    & P^G_{g,t}\in\left[0,P_{\max,g}\right],\ \forall \left(g,t\right)\in G_{\text{th}}\times T,\\[0.1ex]
    &   Q^G_{g,t}\in\left[\min\{Q_{\min,g},0\},\max\{Q_{\max,g},0\}\right],\ \forall \left(g,t\right)\in G_{\text{th}}\times T,\\[0.1ex]
    &P_{g,t}^{\text{res}}\in\left[0,P_{\max,g}- P_{\min,g}\right],\ \forall \left(g,t\right)\in G_{\text{th}}\times T, \\[0.1ex]
    & u_{g,t}, v_{g,t}, w_{g,t} \in \left[0,1\right],\ \forall \left(g,t\right)\in G_{\text{th}}\times T. \\[0.1ex]
\end{array}
\right.
\end{equation}
Also, notice that, when fixing the values of $u_{g,t}, v_{g,t}, w_{g,t}$ in (UC-ACOPF), we obtain a non-convex NLP problem that we denote in this work by ACOPF problem or problem (ACOPF). For clarity, in the remainder of this paper, the decision variables in (RUC-ACOPF) which differ by definition from their counterpart in (UC-ACOPF) are denoted by a tilde; explicitly: $\tilde{u}_{g,t}$, $\tilde{v}_{g,t}$, $\tilde{w}_{g,t}$,  $\tilde{P}^{\text{res}}_{g,t}$, $\tilde{P}^G_{g,t}$  and $\tilde{Q}^G_{g,t}$. 

The natural approach to any relax-and-round strategy can be summarized as follows. First, relax the integrality condition on the binary variables in (UC-ACOPF) to obtain problem (RUC-ACOPF). Next, solve (RUC-ACOPF), and round \textit{somehow} the relaxed variables $\tilde{u}_{g,t}, \tilde{v}_{g,t},\tilde{w}_{g,t}$ to binary values. Then, fix these values in (UC-ACOPF) to obtain problem (ACOPF), and, finally, solve (ACOPF) to attempt to find feasible values for the remaining decision variables.

The key to this strategy is that solving (RUC-ACOPF), a non-convex NLP problem, is conceptually simpler and significantly faster than solving (UC-ACOPF), a non-convex MINLP problem. However, the binary solution obtained after rounding the relaxed solution of (RUC-ACOPF) is typically not suitable for the UC-ACOPF problem. There are two main reasons for this. First, the combinatorial nature of the mixed-integer constraints is lost during relaxation, so rounding alone may not ensure UC feasibility. Second, by relaxing $u_{g,t}\in \{0,1\}$ to $\tilde{u}_{g,t}\in \left[0,1\right]$, $P^G_{g,t}\in\{0\}\times\left[P_{\min,g}, P_{\max,g}\right]$ becomes $\tilde{P}^G_{g,t}\in\left[0, P_{\max,g}\right]$; the domain for the active power transitions from being a unconnected set to becoming a connected one. Here, we propose a slightly different approach. After solving problem (RUC-ACOPF), and previous to rounding, we apply a transformation to $\tilde{u}_{g,t}$ to obtain a more `roundable' variable $\tu_{g,t}^r$ in the sense to be explained in the next section. Additionally, we propose an enhanced rounding formula designed to ensure the feasibility of constraints (\ref{UC_BVL})-(\ref{UC_MTdown}) and better enforce UC-feasibility. However, general feasibility or optimality are not guaranteed; in this sense, we sometimes abuse the notation by referring to the points returned by the algorithms as `solutions'.%
\subsection{Making the relaxed solution `roundable'}
As previously discussed, relaxing $u_{g,t}$ to the interval $\left[0,1\right]$ fundamentally changes the properties of the UC-ACOPF problem. In the original problem, the active power $P^G_{g,t}$ must be either zero or at least $P_{\min,g}$. However, in the relaxed problem (RUC-ACOPF), thermal power plants can inject only the minimum amount of power needed
 to satisfy the demand within the range $\left[0, P_{\max,g}\right]$. Furthermore, while the costs associated with the binary variables $u_{g,t}$, $v_{g,t}$ and $w_{g,t}$ are constant in (\ref{UC_cost}), the relaxed variables $\tilde{u}_{g,t}$, $\tilde{v}_{g,t}$, $\tilde{w}_{g,t}$ $\in\left[0,1\right]$ have linear and progressive costs proportional to their value. When optimizing the relaxed problem, $\tilde{u}_{g,t}$, $\tilde{v}_{g,t}$, $\tilde{w}_{g,t}$ will only approach the value $1$ if this is necessary for feasibility or if it is more cost-effective than remaining closer to zero. The main reason for $\tilde{u}_{g,t}$ tending towards $1$ is to allow $\tilde{P}^G_{g,t}$ to approach $P_{\max,g}$, since
\begin{equation}
    \label{S5eq1}
    P_{\min,g}\tu_{g,t} \le \tP_{g,t}^G \le P_{\max,g}\tu_{g,t},
\end{equation}
assuming $P_{g,t}^{\text{res}}=0$. In conclusion, $\tu_{g,t}$ does not have the same meaning as the binary $u_{g,t}$ since $\tu_{g,t}=1$ does not necessarily mean that the unit is `on', but rather that the plant is producing at maximum capacity. It is not possible to derive the exact on/off status from $\tilde{u}_{g,t}$ as the on/off nature has been relaxed.

{The previous reasons intuitively justify why (RUC-ACOPF) solutions can deviate from the solutions to problem (UC-ACOPF). However, if a direct rounding approach is to be used, we wondered if it might be interesting to redefine \textit{which quantity} from the relaxed solution should be rounded.  If we assume that the main role of $\tu_{g,t}$ is to allow $\tP_{g,t}$ to reach a power level $\tu_{g,t}P_{\max,g}$, and that the relaxed problem retains useful information from the relaxed combinatorial relationships of the original (UC-ACOPF), then we can set
\begin{equation}
    \label{S5eq2_tuMin}
    \tu_{\min,g} = \frac{P_{\min,g}}{P_{\max,g}},
\end{equation}
as a minimum reference for $\tu_{g,t}$, modeling the on/off status of the power plant $g\in G_{\text{th}}$, $\forall t\in T$. If $\tP_{g,t}^G \ge P_{\min,g}$ it is reasonable to expect  $u_{g,t}=1$ in a solution for (UC-ACOPF). For this to happen, taking into account (\ref{S5eq1}), necessarily $\tu_{g,t}\ge P_{min,g}/P_{\max,g}$. Then, if  $\tu_{g,t}\ge  \tu_{\min,g}$, we should round $\tu_{g,t}$ to $1$ when fixing $u_{g,t}$. Although the value (\ref{S5eq2_tuMin}) varies for each thermal power plant, it has the same meaning for all of them; this homogenizes the criteria for rounding since the parameters characterizing each $g\in G_{\text{th}}$ affect the behavior of $\tu_{g,t}$ through the relaxed version of the mixed-integer constraints. }

{Based on the discussion above, we propose two different strategies to modify $\tu_{g,t}$ to obtain better-suited values, $\tu^r_{g,t}$, for rounding to binaries using a direct formula.}

The first strategy, Re-RUC, rescales $\tu_{g,t}$ by the inverse of $\tu_{\min,g}$, yielding
\begin{equation}
    \label{S5eq3_uResc}
    \tu_{g,t}^r = \frac{1}{\tu_{\min,g}}\tu_{g,t} = \frac{P_{\max,g}}{P_{\min,g}}\tu_{g,t}.
\end{equation}
The second strategy, Re-Power, rescales $\tP_{g,t}$ and defines the rounding variable as
\begin{equation}
    \label{S5eq3_uRedf}
    \tu_{g,t}^r = \frac{\tP_{g,t}}{P_{\min,g}}.
\end{equation}
This second approach has been used previously in other works involving rounding procedures. For example in \cite{he2016robust}, to get an initial integer solution for a feasibility pump algorithm. Also in \cite{wang2021optimal}, combined with what we later call `naive rounding formula', rounding values around a fixed threshold.

\begin{figure}[H]
    \centering
    \begin{minipage}{0.33\textwidth}
        \centering
        \includegraphics[width=\textwidth]{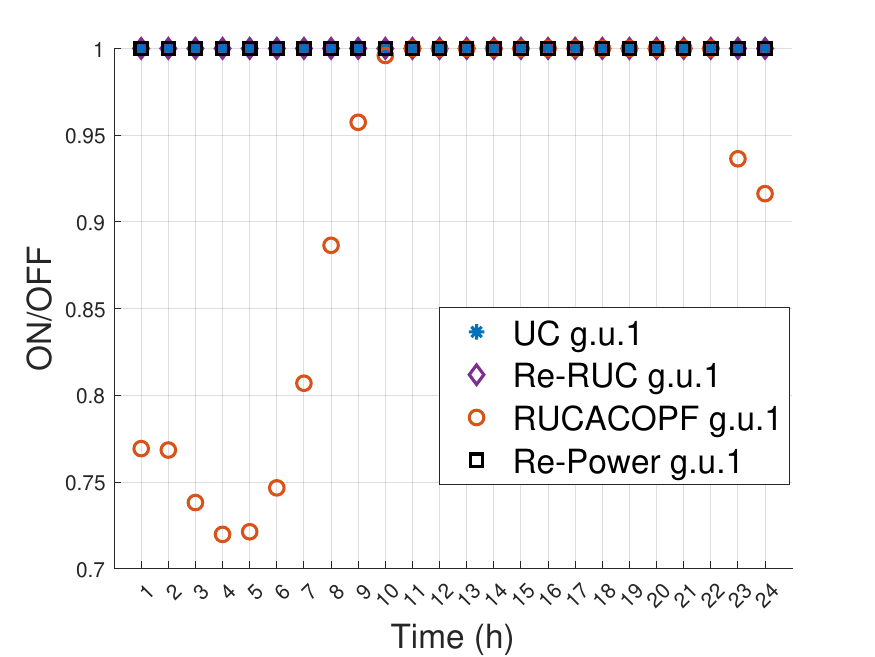}
        \caption*{(a) Generation unit 1.}
    \end{minipage}\hfill
    \begin{minipage}{0.33\textwidth}
        \centering
        \includegraphics[width=\textwidth]{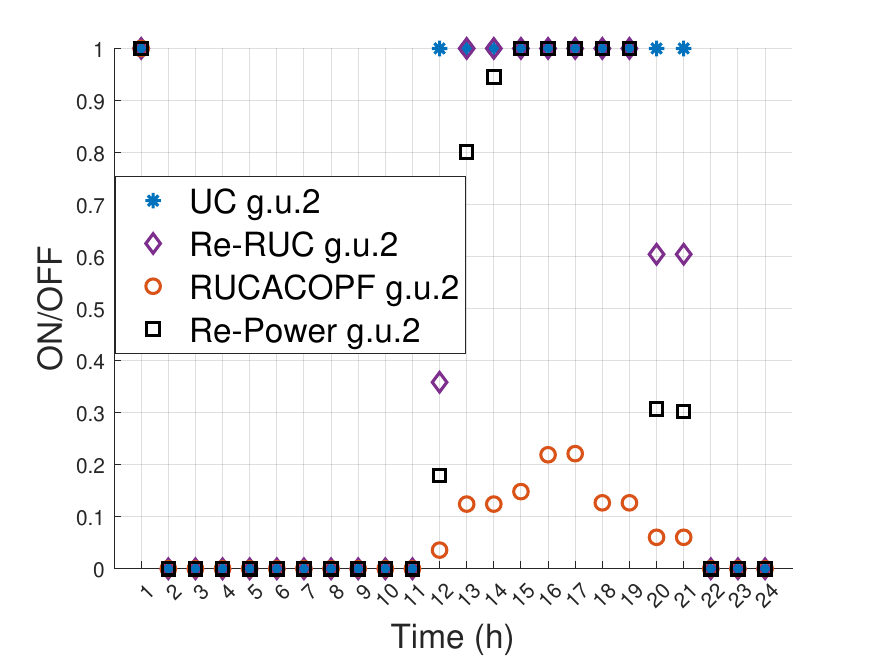}
        \caption*{(b) Generation unit 2.}
    \end{minipage}\hfill
    \begin{minipage}{0.33\textwidth}
        \centering
        \includegraphics[width=\textwidth]{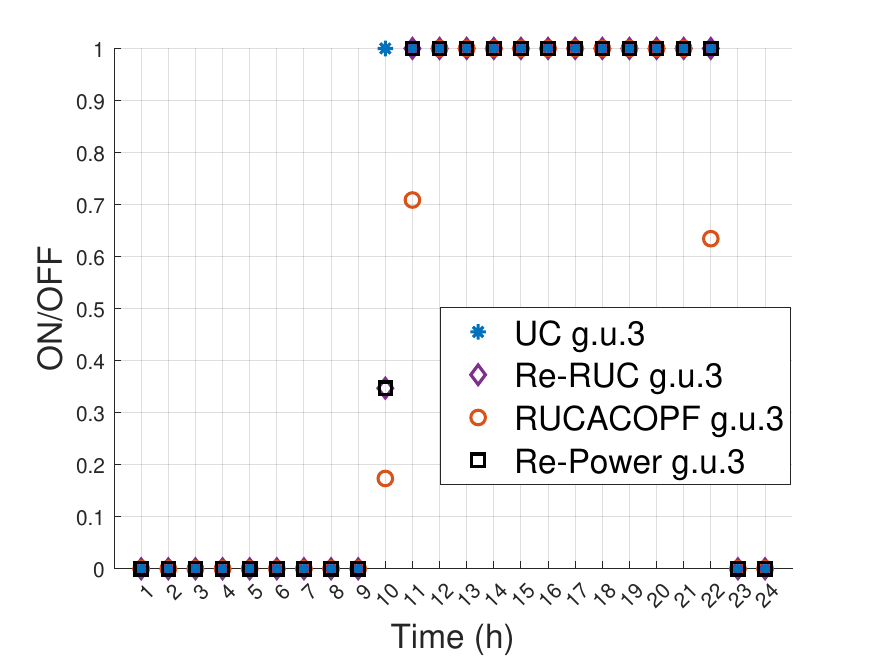}
        \caption*{(c) Generation unit 3.}
    \end{minipage}
    \caption{On/off status of the 3 power plants in the 6-bus test system example.}
    \label{fig:t6b_UCvsRUCvsResc_u}
\end{figure}

We want to numerically evaluate the performance of these strategies by comparing them to the standard approach:
\begin{equation}
    \label{S5eq3_uStan}
    \tu_{g,t}^r = \tu_{g,t}.
\end{equation}
We will detail the numerical experiments in the following sections. To highlight the main features of the strategies, we present some preliminary results related to the rescaling. In Figure \ref{fig:t6b_UCvsRUCvsResc_u}, we show an example for the 6-bus test system from \cite{fu2006ac}, where we compare the integer solution $u_{g,t}$, the relaxed solution $\tu_{g,t}$ obtained using a local NLP solver, and the two proposed rescaling strategies, i.e., Re-RUC and Re-Power. After rescaling, some values of $\tu_{g,t}^r$ may exceed 1. In Figure \ref{fig:t6b_UCvsRUCvsResc_u}, the values $\tu_{g,t}^r > 1$ are rounded down to 1 before plotting. Note that for relaxed values within the range $\left(0,1\right)$, the proposed approaches are expected to make it easier to recover the `desired' integer solution using the natural rounding formula around $0.5$. {As we will see later,  this property of the rescaled relaxed commitments can improve the feasibility rate of solutions returned by heuristics based on relax-and-round procedures. The simplest of these are direct rounding formulas, which we discuss next.}  
%

\subsection{Rounding formulas}

In this section, we will first discuss the standard or naive rounding (NR) formula, as shown in Algorithm \ref{alg:naive_rounding}. The idea behind this approach is simple: when rounding continuous values within the interval $\left(0,1\right)$, a threshold of $0.5$ is used as a reference; values greater than $0.5$ are rounded to $1$, while values less than or equal to $0.5$ are rounded to $0$.

This method poses several challenges when applied to the UC problem since it can easily lead to infeasible solutions. To illustrate these, we refer to the results shown in Figure \ref{fig:t6b_UCvsRUCvsResc_u}, assuming that it is reasonable to aim for recovering a given integer solution. In practice, given a relaxed solution, we will not know a priori to which integer solution we should be trying to round to.
\begin{algorithm}
\small
\caption{Naive rounding formula}\label{alg:naive_rounding}
\begin{algorithmic}[1]
\Require Relaxed and rescaled solution $\tu^r_{g,t}$ and some grid data.
\For{$t \gets 1$ \textbf{to} $t_f$}
    \For{$g \gets 1$ \textbf{to} $N_{\text{g}}$}
        \If{$\tu^r_{g,t} > 0.5$}
            \State $u_{g,t}=1$
        \Else
            \State $u_{g,t}=0$
        \EndIf
    \EndFor
\EndFor
\end{algorithmic}
\end{algorithm}
As noted previously, the NR formula would perform poorly when applied to the relaxed solution $\tu_{g,t}$ (red circle), since some values below $0.5$ would be `erroneously' rounded to zero, even though these units should be committed (turned on). Although our proposed strategies, (\ref{S5eq3_uResc}) and (\ref{S5eq3_uRedf}), produce a more `roundable' relaxed solution $\tu_{g,t}^r$, certain problematic values remain. In particular, the generation of unit 2 at times 12, 20, and 21 and the generation of unit 3 at time 10. In response to such problems, some authors suggest the use of thresholds closer to or even equal to zero to avoid committing units that would not be able to effectively meet the energy demand \cite{kjeldsen2012heuristic}, \cite{sawa2007security}. One way to improve the rounding process, instead of using a fixed threshold such as $0.5$, is to round the relaxed values in descending order or by grouping them into descending intervals (e.g.,$\left[0.9,1\right)$, $\left[0.8,0.9\right)$, etc.). The idea is that processing the higher values first ensures that the most significant generating units are committed, which leads to a more practical and feasible solution compared to using a fixed rounding threshold. However, when using this approach, a stopping criterion is essential to avoid over-committing generating units. A natural criterion would be to stop rounding when the total power generation equals the total demand. The total power demand at time $t\in T$, $P_t^D$, is given by  the sum of the demands at all buses in the system at that, that is,
\begin{equation}
\label{S5:pow_dem}
P_t^D \defeq \sum_{i\in I_{D}} P_{i,t}^D.
\end{equation}
We decided to track the amount of power committed after each rounding decision. When $\tu_{g,t}^r$ is rounded so that $u_{g,t}=1$, thermal generating unit $g\in G_{\text{th}}$ is committed, and its relaxed power generation $\tP_{g,t}^G$ is added to the total committed power $P^G_t$. In other words, the contribution of committed plants at each time step is tracked by adding their power generation from the relaxed solution. This total is then compared with the system demand to ensure sufficient power has been committed, by verifying that $P^G_t \ge P_t^D$. If the committed power meets or exceeds the demand, further rounding is stopped.

Notice that if we apply the Re-Power rescaling (\ref{S5eq3_uRedf}), $\tu_{g,t}^r\ge 1$ ensures that the power generation $\tP_{g,t}^G$ is at least equal to the minimum required power $P_{\min,g}$ for unit $g$, thereby guaranteeing that the unit meets its minimum generation requirement. However, for the Re-RUC strategy (\ref{S5eq3_uResc}), this is not necessarily true, as $\tu_{g,t}$ are involved in additional constraints. Nevertheless, since we are relying on a relax-and-round strategy, it seems reasonable to automatically commit units where $\tu_{g,t}^r\ge 1$.

For the remaining generating units $g \in G_{\text{th}}$ with $\tu_{g,t}^r < 1$ the decision formula must determine whether to round to $1$ (commit the unit) or not. If the unit is committed, it must generate at least its minimum power output, $P_{\min,g}$, according to the constraints of the UC-ACOPF problem. In such a case, rather than simply adding $\tP_{g,t}^G$ to the total committed power $P^G_t$, we should add the greater of $\tP_{g,t}^G$ and $P_{\min,g}$, ensuring that the unit generates enough power to remain within the feasible region:
$$P^G_t= P^G_t+\max\{\tP_{g,t}^G, P_{\min,g}\}.$$
If $P_{\min,g} > \tP_{g,t}^G$, we must also ensure that the additional power, $P_{\min,g} - \tP_{g,t}^G$, is tracked to avoid over-generation. We then want to check whether it is possible to reduce the committed power elsewhere in the system to balance the excess power generated by units operating at their minimum levels. To achieve this, we define a slack power, $P_{s,t}^G$, which initially represents the total power generated above $P_{\min,g}$. When committing a generating unit, the slack power should be updated, for example, by adding $-P_{\min,g} + \tP_{g,t}^G$ (which could be negative). However, to be more accurate, we must account for ramping limits that constrain how much a generating unit's output can change between time steps. Thus, we update the slack power tracker as:
$$ P_{s,t}^G = P_{s,t}^G + \min\{R_g^d-\tP_{g,t-1}^G+\tP_{g,t}^G,- P_{\min,g} + \tP_{g,t}^G\}.$$
Finally, we can track the potential power generation, $P_{p,t}^G$, based on the commitment decisions made up to any point in the rounding process. This quantity can be used to ensure that reserve requirements are met. Assuming a single region for power demand constraints, the condition
$$\sum_{g\in G_{\text{th}}}\tP_{g,t}^{\text{res}}u_{g,t}+P_{p,t}^G > P_{r,t}^{\text{res}}$$
can be used as an additional stopping criterion during the rounding process. After deciding for all $g \in G_{\text{th}}$, $P_{p,t}^G$ can also be used to verify whether the demand can be met, i.e., if $P_{p,t}^G + P_{t}^G > P_t^D$. If this condition is not satisfied, it indicates that the rounding procedure may have failed for that time step.

We have applied all the ideas discussed above to develop a rounding method that we have named the Enhanced Rounding (ER) formula. While it provides a systematic way to round relaxed values, it does not guarantee AC-feasibility for the final solution. Nor does it ensure UC-feasibility, primarily because the combinatorial constraints involving minimum up/down times, $T_{g}^u$ and $T_{g}^d$, are not fully considered. To address the latter, we developed the UC Enhanced Rounding UC-ER formula, outlined in Algorithm \ref{alg:UC-ERF_rounding}. The combinatorial logic is enforced in lines 4-10, while the remaining steps implement the ideas from the ER formula. In this approach, we prioritize rounding by grouping $\tu_{g,t}^r$ into descending intervals, the size of which can be fine-tuned to approximate strict descending order, making the implementation more flexible. The intuitive idea of rounding in descending order, to respect to the relaxed binary values, to match the generation and demand while enforcing minimum up/down times was also mentioned without specifics in \cite{bai2009semi}. It is expected that, once the commitment decisions from the UC-ER formula are fixed and the ACOPF problem is solved, the full set of UC-skeleton and other constraints will be enforced, and a feasible point for the UC-ACOPF problem returned.
\begin{algorithm}[!t]
\caption{UC Enhanced Rounding UC-ER formula}\label{alg:UC-ERF_rounding}
\begin{algorithmic}[1]
\Require Relaxed and rescaled solution $\tu^r_{g,t}$, global power demand for each time $P_t^D$, an array of rounding levels $lvl\left(\cdot\right)$ in $\left[0,1\right]$ of size $N_{lvl}$ and some grid data.
\State Initialize $u_{g,t}=0$ for all $g$ and $t$.
\For{$t \gets 1$ \textbf{to} $t_f$}
    \State Initialize $P_t^G =0 $, $P_{p,t}^G =0 $, $P_{s,t}^G =0$.
    \State Classify the power plants as forced $Cg$ or free $FCg$ to commit.
    \For{$g \in Cg$}
        \State $u_{g,t}=1$,
        \State $P_{s,a}= P_{\min,g}-\tP_{g,t}^G$,
        \State $\tP_{g,t}^G = \max\{P_{\min,g},\tP_{g,t}^G\}$,
        \State $P_t^G = P_t^G + \tP_{g,t}^G$,
        \State $P_{p,t}^G  = P_{p,t}^G + u_{g,t-1}\left(R_g^u-\tP_{g,t}^G -\tilde{P}^{\text{res}}_{g,t}+\tP_{g,t-1}^G\right)
        +\left(1-u_{g,t-1}\right)\left(S_g^u-\tP_{g,t}^G -\tilde{P}^{\text{res}}_{g,t}\right)$,
        \State $P_{s,t}^G  = P_{s,t}^G + \min\{R_g^d-\tP_{g,t-1}^G+\tP_{g,t}^G, P_{s,a}\} $.
    \EndFor
    \State Classify $FCg$ by rounding levels given by $lvl\left(\cdot\right)$.
    \For{$l \gets N_{lvl}$ \textbf{to} $1$}
        \State If {$P_t^G > P_t^D$ and $\sum_{g\in G_{t}} \tP^{\text{res}}_{g,t}u_{g,t} + P_{p,t}^G \ge P^{\text{res}}_{r,t} $}  \textbf{ then $t+1$ and go to }3:
        \For{$g\in lvl\left(l\right)$}
            \State $dP = P_{\min,g} - P_{g,t}$,
            \If{$\left(P_{s,t}^G > dP\right)$ and $\left(P_t^G + P_{\min,g}-(P_{s,t}^G - dP)<PD(t)\right)$}
                \State $u_{g,t}=1$,
                \State $P_{s,a}= P_{\min,g}-\tP_{g,t}^G$,
                \State $\tP_{g,t}^G = \max\{P_{\min,g},\tP_{g,t}\}$,
                \State $P_t^G = P_t^G + \tP_{g,t}^G$,
                \State $P_{p,t}^G  = P_{p,t}^G + u_{g,t-1}\left(R_g^u-\tP_{g,t}^G -\tilde{P}^{\text{res}}_{g,t}+\tP_{g,t-1}^G\right)$
                \State $\quad\quad\quad +\left(1-u_{g,t-1}\right)\left(S_g^u-\tP_{g,t}^G -\tilde{P}^{\text{res}}_{g,t}\right)$,
                \State $P_{s,t}^G  = P_{s,t}^G + \min\{R_g^d-\tP_{g,t-1}^G+\tP_{g,t}^G, P_{s,a}\} $.
            \EndIf
        \EndFor
    \EndFor
\EndFor
\end{algorithmic}
\end{algorithm}
\subsection{Relax-and-round algorithm}
We now formalize our proposed direct approach, leveling on the relaxation and rounding formulas discussed in the previous sections. The main steps of the methodology are summarized in Algorithm \ref{alg:relax_round}. In \text{Step 2}, one can choose between: the standard approach (\ref{S5eq3_uStan}), the rescaling for $\tu_{g,t}$ (\ref{S5eq3_uResc}) or the rescaling based on $\tP_{g,t}^G$ (\ref{S5eq3_uRedf}). In \text{Step 3} the NR, ER, or UC-ER formulas can be applied. In general, any other appropriate rounding approach can be used.

\begin{algorithm}[H]
\small
\caption{Relax-and-round with formulas}\label{alg:relax_round}
\begin{algorithmic}[1]
\State Consider (UC-ACOPF).
\State \text{Step 1}: Solve (RUC-ACOPF) obtainning $\tilde{u}_{g,t}$.
\State \text{Step 2}: Rescale $\tilde{u}_{g,t}$ to obtain $\tu^r_{g,t}$.
\State \text{Step 3}: Round $\tu^r_{g,t}$ using a rounding formula. Let $\boldsymbol{y}_{\text{R}}$ be then the integer solution.
\State{\text{Step 4}: Fix the integer variables $\boldsymbol{y}=\boldsymbol{y}_{\text{R}}$ in (UC-ACOPF), obtaining (ACOPF). Solve (ACOPF). Let $\boldsymbol{x}_{\text{AC}}$ be its solution.}
\State{Return $\left(\boldsymbol{x}_{\text{AC}},\boldsymbol{y}_{\text{R}}\right)$.}
\end{algorithmic}
\end{algorithm}
\section{Numerical results} \label{Sect3}
{This section presents numerical results for the (UC-ACOPF) problem, applying the direct relax-and-round approach introduced in the previous section to two common test systems from the literature: the 6-bus test system and a modified IEEE 118-bus system. The goal is to assess the feasibility of the solutions returned by the algorithms when using the proposed rescaling strategies and rounding formulas.}

To solve the (RUC-ACOPF) and (ACOPF) problems, we employ classic penalized sequential linear programming (PSLP) algorithms (see, for example, \cite{bazaraa2006nonlinear}, from which we adopted the parameter notation). In a PSLP algorithm, the nonlinear constraints are linearized using first-order Taylor expansion. To ensure feasibility throughout the iteration process, slack variables are incorporated into the linearized constraints of non-convex equations. These slack variables are then penalized in the cost function. The initial guess for the decision variables in \text{Step 1} of Algorithm \ref{alg:relax_round} is a modified flat value for the complex voltage and the generation commitment at $t_0$ for the power plants across the 24-hour time frame. In particular, we set $V_{i,t}=\left(V_{\max,i}+V_{\min,i}\right)/2$, for all $t\in T$; this change helps prevent numerical issues when initializing the SLP algorithm, which can occur if the voltage is set to $1$ {in per-unit system (p.u.)}
and this value is too close to the maximum or minimum allowed voltage limits for any bus $i$. In \text{Step 4}, we initialize the algorithm with a modified flat value for the complex voltage and set the power generation commitments based on the binary values obtained from the rounding formula. Once these binary decisions are fixed, we assign appropriate initial values to the active and reactive power for each committed power plant. For example, the reactive power $Q_{t,g}^G$ is set to $0$ for all power plants, and active power $P_{t,g}^G$ is set close to the plant’s minimum operating level $P_{\min,g}$ when the plant is committed (i.e., $u_{g,t}=1$) or $0$ if the plant is not committed (i.e., $u_{g,t}=0$). This ensures that the basic UC constraints are satisfied. Despite this setup can be AC-infeasible, the algorithm can still proceed because it introduces penalized slack variables, which allow flexibility in meeting these AC constraints during the iterations. The PSLP method is a local minimizer, with several parameters involved that affect both the performance of the algorithm and the solutions obtained. The parameter values have been summarized in Table \ref{tab:PSLP_parameters}. Additionally, we have used various penalty weights $\{5.\mathrm{e}{3},5.\mathrm{e}{4},5.\mathrm{e}{5},5.\mathrm{e}{6}, 5.\mathrm{e}{7}\}$ in \text{Step 1} of the PSLP to solve the (RUC-ACOPF) problem, which resulted in different relaxed values of $\tu_{g,t}$ and, consequently, different integer solutions after rounding. These penalty weight values were selected based on empirical numerical experimentation.

\begin{table}[ht]
    \centering
    \begin{tabular}{|c|c|c|c|c|c|c|}
        \hline
        Stop $\epsilon$ & Penalty & $\rho_0$ & $\rho_1$ & $\rho_2$ & $\beta$ & $\Delta_{LB}$ \\ \hline
        $1.\mathrm{e}{-4}$ & $5.\mathrm{e}{6}$ & $1.\mathrm{e}{-6}$ & $0.25$ & $0.75$ & $0.5$ & $1.\mathrm{e}{-6}$ \\ \hline
    \end{tabular}
    \caption{Values of the parameters for the PSLP algorithms.}
    \label{tab:PSLP_parameters}
\end{table}

All the linear problems involved are solved using Gurobi v9.5.2 \cite{gurobi}. The feasibility of the solutions is double-checked using SCIP \cite{Achterberg2009} through NEOS Server \cite{czyzyk_et_al_1998} or MAinGO \cite{najman2019mccormick}. For all solvers, the standard configuration is used. The feasibility tolerance is set to $1.\mathrm{e}{-6}$. The PSLP parameters were chosen to align with this threshold, and we consider solutions to be (numerically) AC-feasible when constraints (\ref{ACPF_PV_APB})-(\ref{ACPF_PV_RPB}) and (\ref{ACPF_PV_LineLim}) are satisfied within this tolerance. As usual, the ACPF equations being solved are first normalized using the per-unit system.

\subsection{6-bus test system}
\begin{figure}[th]
 \centering
 \includegraphics[scale=0.4]{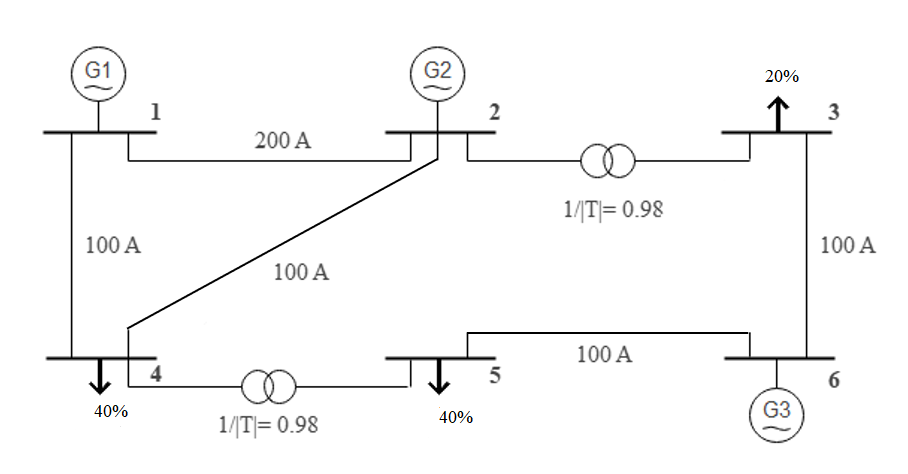}
  \caption{6-bus test system circuit.}\label{fig:6busgraph}
\end{figure}

The first example we present is the 6-bus test system from \cite{fu2005security}. The system circuit modeling the power grid is represented in Figure \ref{fig:6busgraph}. It consists of 6 buses and 7 elements: 2 tap-changing transformers and 5 transmission lines. There are 3 thermal generators placed at buses 1, 2, and 6, and loads or power demands on buses 3, 4, and 5. Despite its small size, 216 binary and 1344 continuous variables, global solvers available in NEOS and MAiNGO fail to solve it within 8 hours using their standard configurations..

Table \ref{tab:6bus_summary} summarizes the performance of our direct strategy, using different combinations of rescaling approaches and rounding formulas. Results are presented only for a penalty  weight of $5.\mathrm{e}{6}$ in the PSLP used to solve problem (RUC-ACOPF), as similar outcomes were observed for other values. In particular, the NR formula was unable to produce a feasible solution within the specified tolerance of $1.\mathrm{e}{-6}$, though the Re-RUC approach showed better feasibility compared to the other two methods. For the ER and UC-ER rounding strategies, the rescaling approach had no significant impact on the final result, as both consistently produced the same feasible UC-ACOPF solution.  The integer commitments returned by the UC-ER and ER formulas are detailed in Table \ref{tab:6bus_solcomp} and compared to UC-ACOPF and UC with DC power flow constraints (UC-DCOPF) solutions for the same test system from the literature. The power plant commitment and ACPF state costs are consitent with the global solution provided in \cite{liu2018global}; moreover, they are comparable to solutions obtained using local minimizers as in \cite{castillo2016unit}, allowing us to compare the UC solution and the relaxed solution of the RUC. In this case, the binary variables and their relaxed counterparts are directly linked through the UC-ER and ER rounding formulas, i.e., applying these formulas to the relaxed solutions yields this same integer solution.

\begin{figure}[ht]
    \centering
    \begin{subfigure}{0.5\textwidth}
        \centering
        \includegraphics[width=\textwidth]{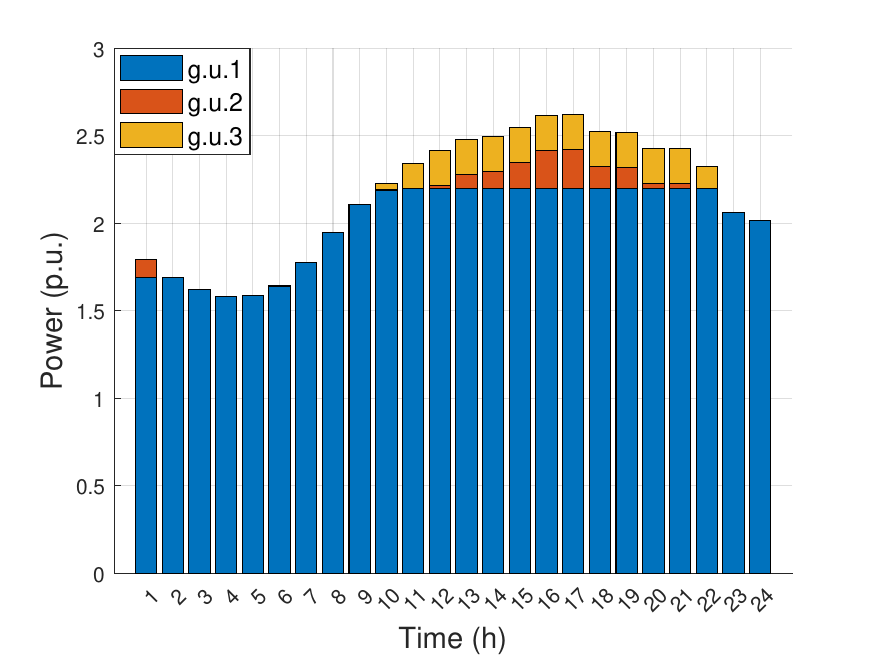}
        \caption{(RUC-ACOPF)}
    \end{subfigure}\hfill
    \begin{subfigure}{0.5\textwidth}
        \centering
        \includegraphics[width=\textwidth]{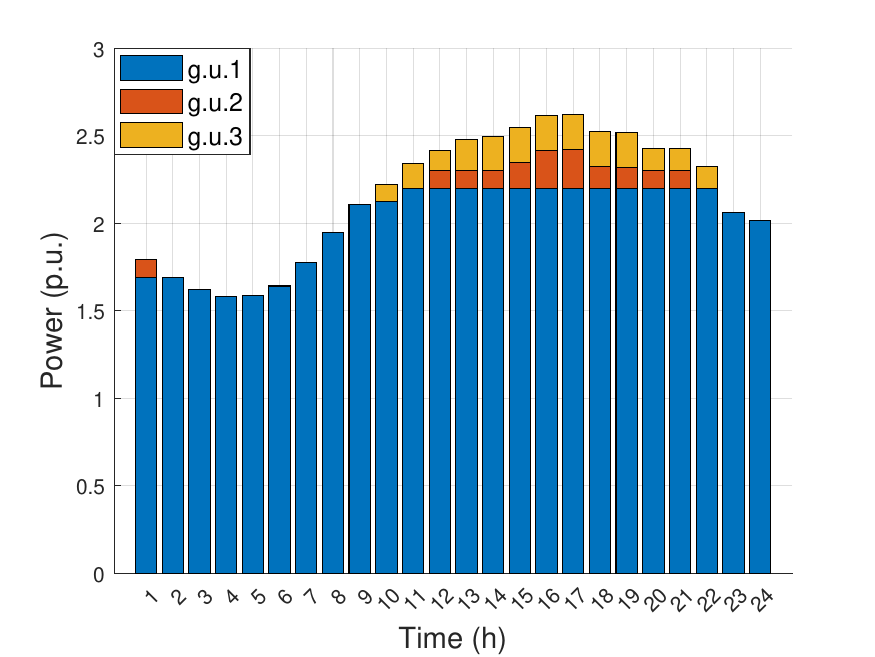}
        \caption{(UC-ACOPF)}
    \end{subfigure}
    \caption{Active power generation in the 6-bus test system example.}\label{fig:6bus_Pg}
\end{figure}

\begin{figure}[ht]
    \begin{subfigure}{0.5\textwidth}
        \centering
        \includegraphics[width=\textwidth]{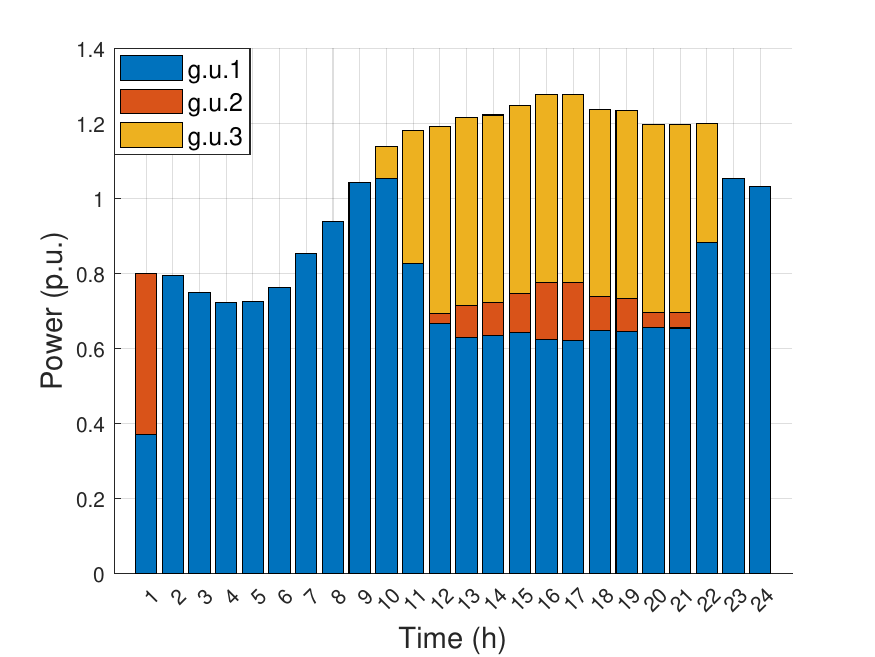}
        \caption{(RUC-ACOPF)}
    \end{subfigure}\hfill
    \begin{subfigure}{0.5\textwidth}
        \centering
        \includegraphics[width=\textwidth]{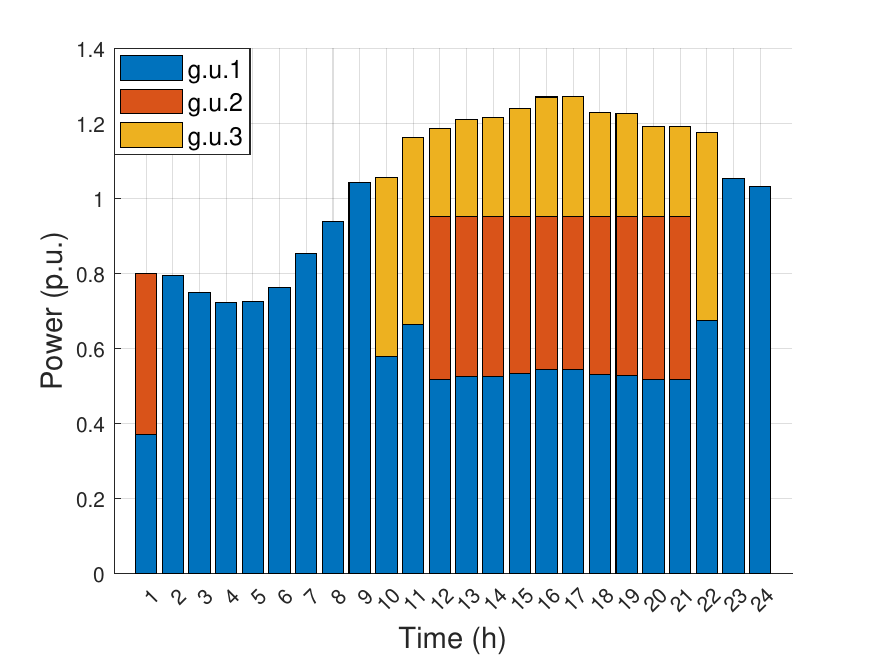}
        \caption{(UC-ACOPF)}
    \end{subfigure}
   \caption{Reactive power generation in the 6-bus test system example.}\label{fig:6bus_Qg}
\end{figure}

\begin{table}[!t]
\centering
\begin{tabular}{lccccccccc}
\toprule
\multicolumn{2}{c}{}&  \multicolumn{3}{c}{$1^{\textrm{st}}$ PSLP }
     & \multicolumn{2}{c}{$2^{\textrm{nd}}$ PSLP } &  \multicolumn{3}{c}{Results} \\
    \midrule
  \textbf{Rescale} &  \textbf{Round} &  Penalty & Iter. &  AC-feas. & Iter. &  AC-feas. & Feas. & Cost($\$$) & Time (s)  \\
\toprule
None & Naive & $5.\mathrm{e}{6}$ & 5 & $3.16\mathrm{e}{-13}$ & 17 & $4.12\mathrm{e}{1}$ & \xmark & $94704$ & $2.7$ \\
Re-RUC & Naive & $5.\mathrm{e}{6}$ & 5 & $3.16\mathrm{e}{-13}$ & 32 & $6.17\mathrm{e}{-2}$ & \xmark & $101058$ & $3.8$ \\
Re-Power & Naive & $5.\mathrm{e}{6}$ & 5 & $3.16\mathrm{e}{-13}$ & 30 & $7.53\mathrm{e}{1}$ & \xmark & $100230$ & $3.6$ \\ \hline
None & ER & $5.\mathrm{e}{6}$ & 5 & $3.16\mathrm{e}{-13}$ & 34 & $1.86\mathrm{e}{-8}$ & \cmark & $101763$ & $4.0$ \\
Re-RUC & ER & $5.\mathrm{e}{6}$ & 5 & $3.16\mathrm{e}{-13}$ & 34 & $1.86\mathrm{e}{-8}$ & \cmark & $101763$ & $4.0$ \\
Re-Power & ER & $5.\mathrm{e}{6}$ & 5 & $3.16\mathrm{e}{-13}$ & 34 & $1.86\mathrm{e}{-8}$ & \cmark & $101763$ & $3.9$ \\ \hline
None & UC-ER & $5.\mathrm{e}{6}$ & 5 & $3.16\mathrm{e}{-13}$ & 34 & $1.86\mathrm{e}{-8}$ & \cmark & $101763$ & $4.0$ \\
Re-RUC & UC-ER & $5.\mathrm{e}{6}$ & 5 & $3.16\mathrm{e}{-13}$ & 34 & $1.86\mathrm{e}{-8}$ & \cmark & $101763$ & $4.0$ \\
Re-Power & UC-ER & $5.\mathrm{e}{6}$ & 5 & $3.16\mathrm{e}{-13}$ & 34 & $1.86\mathrm{e}{-8}$ & \cmark & $101763$ & $4.0$ \\
\bottomrule
\end{tabular}
\caption{Summary of the direct approaches performance for the 6-bus test system example.}
\label{tab:6bus_summary}
\end{table}
\begin{table}
\centering
\begin{tabular}{cccccc}
\hline  & \multicolumn{2}{c}{\cite{castillo2016unit}} & \cite{liu2018global}  & \multicolumn{2}{c}{Rescale + Round strategy} \\

  & (UC) & (UC-DCOPF) & (UC-ACOPF) & Re-Power+UC-ER$^{(*)}$& Re-RUC+NR \\
\hline
 Cost($\$$) & 101269 & 106887 & 101763 & 101763 \cmark &  101058  \xmark \\
g.u. 1 & 1-24 & 1-24 & 1-24 & 1-24 &  1-24  \\
g.u. 2 & 1,12-21 & 1,11-22 & 1,12-21 & 1,12-21 & 1,13-21 \\
g.u. 3 & 10-22 & 10-22 & 10-22 & 10-22 & 11-22 \\
 \hline
  \multicolumn{6}{l}{\scriptsize (*) Also, all combinations including the different rescale and the ER and UC-ER formulas.} \\
\end{tabular}
\caption{Comparison between commitment and cost solutions for the 6-bus test system example.}
\label{tab:6bus_solcomp}
\end{table}
{ The effects of the relaxation on the binary variables for the 6-bus test system was previously illustrated in Figure \ref{fig:t6b_UCvsRUCvsResc_u}.  In the following, we compare the decision variables that characterize the AC power flow state, the power balance at the buses, and the power flow through the lines for the local solutions obtained for (RUC-ACOPF) and (UC-ACOPF), intending to highlight the most significant qualitative differences between the two solutions.} Given that both problems are non-convex, a direct comparison of their AC state solutions is generally not very informative, as non-convexity often leads to multiple local solutions. Nevertheless, the RUC problem produces smoother profiles over time than the UC problem. This is because the integrality relaxation in the RUC softens the discrete jumps in the binary decision variables, particularly those governing the generator on/off states, which in turn smooths the time-coupling effects between decision variables. Consequently, while the UC problem exhibits abrupt changes in generation and power flow due to the binary nature of its constraints, the RUC problem shows more gradual transitions over time. Figure \ref{fig:6bus_Pg} illustrates this by comparing both solutions' total active power generation. While the generation profiles are similar, the RUC solution does not enforce the minimum generation constraints found in the UC problem.  Instead, the RUC solution allows generating units to produce less than their minimum output, prioritizing the cheaper units (in ascending order: g.u. 1, 3, and 2 considering raw generation cost, assuming that power plant 1 is already on). Any remaining power demand is then met by committing the more expensive generating units, creating a difference in how active power is distributed across units for the two problems. In Figure \ref{fig:6bus_Qg}, we can see that the overall reactive power generation profile differs significantly between the two solutions. This discrepancy affects the power balance in the grid and leads to qualitative and quantitative differences in bus voltages and power flow through the lines, which explains the difficulties in achieving AC-feasible solutions when solving problem (ACOPF) after rounding the relaxed unit commitment.

\subsection{118-bus test system}
In this section, we present numerical results for the modified IEEE 118-bus test system from \cite{castillo2016unit}. It consists of 118 buses, 54 thermal generators, 14 synchronous condensers, and 186 elements: 8 tap-changing transformers, 1 phase shifter, and 177 transmission lines. This medium-size benchmark is widely used to evaluate the scalability of solvers for the UC-ACOPF problem. No MINLP solver available on the NEOS server was able to return a solution after 8 hours of computation for this test using their standard configurations, which is consistent with the results reported in \cite{tejada2019unit}. First, we provide results where the descending and rounding procedure on $\tu_{g,t}^r$ was applied based on strict value ordering. Moreover, the results were shown for all penalty weight values specified at the beginning of Section \ref{Sect3}.

\begin{table}[!t]
\centering
\begin{tabular}{lcccccccc}
\toprule
\multicolumn{9}{c}{ \textbf{Naive Round (NR)}} \\
\toprule
\multicolumn{1}{c}{}&  \multicolumn{3}{c}{$1^{\textrm{st}}$ PSLP }
     & \multicolumn{2}{c}{$2^{\textrm{nd}}$ PSLP } &  \multicolumn{3}{c}{Results} \\
    \midrule
  \textbf{Rescale}  &  Penalty & Iter. &  AC-feas. & Iter. &  AC-feas. & Feas. & Cost($\$$) & Time (s)  \\
\toprule
None  & $5.\mathrm{e}{3}$ & 276 & $9.13\mathrm{e}{-6}$ & 123 & $1.24\mathrm{e}{+1}$ & \xmark & $740619$ & $993$ \\
Re-RUC & $5.\mathrm{e}{3}$ & 276 & $9.13\mathrm{e}{-6}$ & 40 & $2.71\mathrm{e}{-8}$ & \cmark & $847460$ & $850$ \\
Re-Power & $5.\mathrm{e}{3}$ & 276 & $9.13\mathrm{e}{-6}$ & 41& $2.61\mathrm{e}{-8}$ & \cmark & $846272$ & $852$ \\
\hline
None & $5.\mathrm{e}{4}$ & 40 & $5.51\mathrm{e}{-6}$ & 351 & $1.07\mathrm{e}{+1}$ & \xmark & $741383$ & $717$ \\
Re-RUC & $5.\mathrm{e}{4}$ & 40 & $5.51\mathrm{e}{-6}$ & 41 & $2.65\mathrm{e}{-8}$ & \cmark & $847633$ & $197$ \\
Re-Power & $5.\mathrm{e}{4}$ & 40 & $5.51\mathrm{e}{-6}$ & 41 & $2.65\mathrm{e}{-8}$ & \xmark & $846332$ & $192$ \\
\hline
None & $5.\mathrm{e}{5}$ & 36 & $1.10\mathrm{e}{-6}$ & 370 & $1.01\mathrm{e}{+1}$ & \xmark & $760409$ & $761$ \\
Re-RUC & $5.\mathrm{e}{5}$ & 36 & $1.10\mathrm{e}{-6}$ & 39 & $2.64\mathrm{e}{-8}$ & \cmark & $847742$ & $179$ \\
Re-Power & $5.\mathrm{e}{5}$ & 36 & $1.10\mathrm{e}{-6}$ & 41 & $2.61\mathrm{e}{-8}$ & \xmark & $846432$ & $184$ \\
\hline

None & $5.\mathrm{e}{6}$ & 39 & $2.46\mathrm{e}{-8}$ & 184 & $1.02\mathrm{e}{+1}$ & \xmark & $759607$ & $438$ \\
Re-RUC & $5.\mathrm{e}{6}$ & 39 & $2.46\mathrm{e}{-8}$ & 40 & $2.77\mathrm{e}{-8}$ & \cmark & $847740$ & $191$ \\
Re-Power & $5.\mathrm{e}{6}$ & 39 & $2.46\mathrm{e}{-8}$ & 41 & $2.58\mathrm{e}{-8}$ & \cmark & $846543$ & $194$ \\
\hline
None & $5.\mathrm{e}{7}$ & 44 & $3.45\mathrm{e}{-9}$ & 228 & $9.97$ & \xmark & $740015$ & $529$ \\
Re-RUC & $5.\mathrm{e}{7}$ & 44 & $3.45\mathrm{e}{-9}$ & 41 & $2.56\mathrm{e}{-8}$ & \cmark & $847713$ & $209$ \\
Re-Power & $5.\mathrm{e}{7}$ & 44 & $3.45\mathrm{e}{-9}$ & 41 & $2.67\mathrm{e}{-8}$ & \xmark & $846527$ & $211$ \\
\bottomrule
\end{tabular}
    \caption{Summary of the NR formula performance for the 118-bus test system example.}
    \label{tab:118bus_summary_naive}
\end{table}

Table \ref{tab:118bus_summary_naive} summarizes the performance of the NR formula. The rescaling approach had a clear impact on the solution feasibility. Without rescaling $\tu_{g,t}^r$, no feasible solution was obtained. With the Re-RUC approach, all returned integer solutions were UC-feasible, and an AC-feasible state was found when solving the (ACOPF) problem in \text{Step 4}. Conversely, the Re-Power approach achieved feasible points for some penalty weights, but not for all. The solution costs were consistently lower when using Re-Power compared to Re-RUC. In summary, Re-RUC rescaling proves to be the best option for the NR formula in terms of consistently producing feasible solutions.
\begin{table}[!t]
\centering
\begin{tabular}{lcccccccc}
\toprule
\multicolumn{9}{c}{ \textbf{Enhanced Round (ER)}} \\
\toprule
\multicolumn{1}{c}{}&  \multicolumn{3}{c}{$1^{\textrm{st}}$ PSLP }
     & \multicolumn{2}{c}{$2^{\textrm{nd}}$ PSLP } &  \multicolumn{3}{c}{Results} \\
    \midrule
  \textbf{Rescale}  &  Penalty & Iter. &  AC-feas. & Iter. &  AC-feas. & Feas. & Cost($\$$) & Time (s)  \\
\toprule
None & $5.\mathrm{e}{3}$ & 276 & $9.13\mathrm{e}{-6}$ & 104 & $2.71\mathrm{e}{-8}$ & \xmark & $852542$ & $942$ \\
Re-RUC & $5.\mathrm{e}{3}$ & 276 & $9.13\mathrm{e}{-6}$ & 39 & $2.86\mathrm{e}{-8}$ & \xmark & $854045$ & $836$ \\
Re-Power & $5.\mathrm{e}{3}$ & 276 & $9.13\mathrm{e}{-6}$ & 41 & $2.68\mathrm{e}{-8}$ & \xmark & $852001$ & $840$ \\
\hline
None & $5.\mathrm{e}{4}$ & 40 & $5.51\mathrm{e}{-6}$ & 41 & $2.77\mathrm{e}{-8}$ & \xmark & $853636$ & $190$ \\
Re-RUC & $5.\mathrm{e}{4}$ & 40 & $5.51\mathrm{e}{-6}$ & 41 & $2.80\mathrm{e}{-8}$ & \xmark & $853120$ & $192$ \\
Re-Power & $5.\mathrm{e}{4}$ & 40 & $5.51\mathrm{e}{-6}$ & 42 & $2.41\mathrm{e}{-8}$ & \xmark & $851622$ & $194$ \\
\hline
None & $5.\mathrm{e}{5}$ & 36 & $1.10\mathrm{e}{-6}$ & 37 & $2.67\mathrm{e}{-8}$ & \xmark & $852722$ & $174$ \\
Re-RUC & $5.\mathrm{e}{5}$ & 36 & $1.10\mathrm{e}{-6}$ & 41 & $2.46\mathrm{e}{-8}$ & \xmark & $852920$ & $179$  \\
Re-Power & $5.\mathrm{e}{5}$ & 36 & $1.10\mathrm{e}{-6}$ & 41 & $2.49\mathrm{e}{-8}$ & \xmark & $851115$ & $181$ \\
\hline
None & $5.\mathrm{e}{6}$ & 39 & $2.46\mathrm{e}{-8}$ & 42 & $2.86\mathrm{e}{-8}$ & \xmark & $854296$ & $193$ \\
Re-RUC & $5.\mathrm{e}{6}$ & 39 & $2.46\mathrm{e}{-8}$ & 40 & $2.56\mathrm{e}{-8}$ & \xmark & $852668$ & $191$ \\
Re-Power & $5.\mathrm{e}{6}$ & 39 & $2.46\mathrm{e}{-8}$ & 42 & $2.39\mathrm{e}{-8}$ & \xmark & $851172$ & $194$ \\
\hline
None & $5.\mathrm{e}{7}$ & 44 & $3.45\mathrm{e}{-9}$ & 40 & $2.6\mathrm{e}{-8}$ & \xmark & $853136$ & $206$ \\
Re-RUC & $5.\mathrm{e}{7}$ & 44 & $3.45\mathrm{e}{-9}$ & 43 & $2.76\mathrm{e}{-8}$ & \xmark & $855147$ & $211$ \\
Re-Power & $5.\mathrm{e}{7}$ & 44 & $3.45\mathrm{e}{-9}$ & 40 & $2.48\mathrm{e}{-8}$ & \xmark & $851907$ & $208$ \\
\bottomrule
\end{tabular}
    \caption{Summary of the ER formula performance for the 118-bus test system example.}
    \label{tab:118bus_summary_ER}
\end{table}
Table \ref{tab:118bus_summary_ER} shows the results for the ER formula. Although all the solutions returned were AC-feasible, none were UC-feasible. Note that in order to recover UC-feasible solutions using the NR formula, it was sufficient to apply the Re-RUC approach. This suggests that the ER formula may be overwriting some of the useful intrinsic information in the relaxed and rescaled solution.
\begin{table}[!t]
\centering
\begin{tabular}{lcccccccc}
\toprule
\multicolumn{9}{c}{ \textbf{UC Enhanced Round (UC-ER)}} \\
\toprule
\multicolumn{1}{c}{}&  \multicolumn{3}{c}{$1^{\textrm{st}}$ PSLP }
     & \multicolumn{2}{c}{$2^{\textrm{nd}}$ PSLP } &  \multicolumn{3}{c}{Results} \\
    \midrule
  \textbf{Rescale}  &  Penalty & Iter. &  AC-feas. & Iter. &  AC-feas. & Feas. & Cost($\$$) & Time (s)  \\
\toprule
None & $5.\mathrm{e}{3}$ & 276 & $9.13\mathrm{e}{-6}$ & 42 & $2.52\mathrm{e}{-8}$ & \cmark & $850421$ & $833$ \\
Re-RUC & $5.\mathrm{e}{3}$ & 276 & $9.13\mathrm{e}{-6}$ & 40 & $2.6\mathrm{e}{-8}$ & \cmark & $851634$ & $835$ \\
Re-Power & $5.\mathrm{e}{3}$ & 276 & $9.13\mathrm{e}{-6}$ & 41 & $2.54\mathrm{e}{-8}$ & \cmark & $849031$ & $837$ \\
\hline
None & $5.\mathrm{e}{4}$ & 40 & $5.51\mathrm{e}{-6}$ & 38 & $2.48\mathrm{e}{-8}$ & \cmark & $851200$ & $187$ \\
Re-RUC & $5.\mathrm{e}{4}$ & 40 & $5.51\mathrm{e}{-6}$ & 40 & $2.54\mathrm{e}{-8}$ & \cmark & $851150$ & $187$ \\
Re-Power & $5.\mathrm{e}{4}$ & 40 & $5.51\mathrm{e}{-6}$ & 39 & $2.28\mathrm{e}{-8}$ & \cmark & $848672$ & $188$ \\
\hline
None & $5.\mathrm{e}{5}$ & 36 & $1.10\mathrm{e}{-6}$ & 39 & $2.35\mathrm{e}{-8}$ & \cmark & $850835$ & $179$ \\
Re-RUC & $5.\mathrm{e}{5}$ & 36 & $1.10\mathrm{e}{-6}$ & 41 & $2.70\mathrm{e}{-8}$ & \cmark & $851328$ & $181$ \\
Re-Power & $5.\mathrm{e}{5}$ & 36 & $1.10\mathrm{e}{-6}$ & 39 & $2.32\mathrm{e}{-8}$ & \cmark & $849138$ & $178$ \\
\hline
None & $5.\mathrm{e}{6}$ & 39 & $2.46\mathrm{e}{-8}$ & 40 & $2.59\mathrm{e}{-8}$ & \cmark & $850817$ & $190$ \\
Re-RUC & $5.\mathrm{e}{6}$ & 39 & $2.46\mathrm{e}{-8}$ & 38 & $2.54\mathrm{e}{-8}$ & \cmark & $851456$ & $188$ \\
Re-Power & $5.\mathrm{e}{6}$ & 39 & $2.46\mathrm{e}{-8}$ & 41 & $2.61\mathrm{e}{-8}$ & \cmark & $848837$ & $193$ \\
\hline
None & $5.\mathrm{e}{7}$ & 44 & $3.45\mathrm{e}{-9}$ & 38 & $2.38\mathrm{e}{-8}$ & \cmark & $850980$ & $205$ \\
Re-RUC & $5.\mathrm{e}{7}$ & 44 & $3.45\mathrm{e}{-9}$ & 40 & $2.61\mathrm{e}{-8}$ & \cmark & $851463$ & $208$ \\
Re-Power & $5.\mathrm{e}{7}$ & 44 & $3.45\mathrm{e}{-9}$ & 39 & $2.52\mathrm{e}{-8}$ & \cmark & $849437$ & $207$ \\
\bottomrule
\end{tabular}
    \caption{Summary of the UC-ER formula performance for the 118-bus test system example.}
    \label{tab:118bus_summary_UCER}
\end{table}
Table \ref{tab:118bus_summary_UCER} summarizes the results for the UC-ER formula. All solutions obtained were feasible. The rescaling approach did not affect the feasibility, but the Re-Power method systematically yielded the least costly feasible solution.
\begin{table}[!t]
\centering
\begin{tabular}{lcccccccc}
\toprule
\multicolumn{9}{c}{ \textbf{UC Enhanced Round (UC-ER)}} \\
\toprule
\multicolumn{1}{c}{} & \multicolumn{1}{c}{ } & \multicolumn{2}{c}{$1^{\textrm{st}}$ PSLP }
     & \multicolumn{2}{c}{$2^{\textrm{nd}}$ PSLP } &  \multicolumn{3}{c}{Results} \\
    \midrule
  \textbf{Rescale}  &  Inter. width & Iter. &  AC-feas. & Iter. &  AC-feas. & Feas. & Cost($\$$) & Time (s)  \\
\toprule
None & $0.05$ & 39 & $2.46\mathrm{e}{-8}$ & 38 & $2.43\mathrm{e}{-8}$ & \cmark & $850938$ & $189$ \\
Re-RUC & $0.05$ & 39 & $2.46\mathrm{e}{-8}$ & 39 & $2.53\mathrm{e}{-8}$ & \cmark & $851503$ & $191$ \\
Re-Power & $0.05$ & 39 & $2.46\mathrm{e}{-8}$ & 41 & $2.61\mathrm{e}{-8}$ & \cmark & $848837$ & $193$ \\
\hline
None & $0.1$ & 39 & $2.46\mathrm{e}{-8}$ & 40 & $2.49\mathrm{e}{-8}$ & \cmark & $852272$ & $190$ \\
Re-RUC & $0.1$ & 39 & $2.46\mathrm{e}{-8}$ & 39 & $2.53\mathrm{e}{-8}$ & \cmark & $851503$ & $187$ \\
Re-Power & $0.1$ & 39 & $2.46\mathrm{e}{-8}$ & 41 & $2.61\mathrm{e}{-8}$ & \cmark & $848837$ & $193$ \\
\hline
None & $0.2$ & 39 & $2.46\mathrm{e}{-8}$ & 38 & $7.3\mathrm{e}{-3}$ & \xmark & $858559$ & $189$ \\
Re-RUC & $0.2$ & 39 & $2.46\mathrm{e}{-8}$ & 38 & $2.53\mathrm{e}{-8}$ & \cmark & $851772$ & $185$ \\
Re-Power & $0.2$ & 39 & $2.46\mathrm{e}{-8}$ & 41 & $2.56\mathrm{e}{-8}$ & \cmark & $849525$ & $193$ \\
\bottomrule
\end{tabular}
    \caption{Performance of UC-ER formula on the 118-bus test system based on varying grouping interval widths.}
    \label{tab:118bus_summary_UCER2}
\end{table}

Next, we set the penalty weight for the PSLP in \text{Step 1} to $5.\mathrm{e}{6}$ and repeated the experiments, this time by dividing the $\tu_{g,t}^r$ values into different intervals, with each subsequent experiment increasing the interval size. The order of each generator within these intervals was determined by the grid numbering system, ensuring a structured grouping. The results of these experiments are detailed in Table \ref{tab:118bus_summary_UCER2}. We observed that as the interval size increased, the overall robustness of the process decreased. However, since rounding formulas are not optimization approaches, rounding in descending order based on the relaxed commitment value can lead to low-quality, suboptimal solutions. A grouping of the rescaled values with reasonable interval widths in combination with a priority list approach could improve the quality of the returned solutions in terms of cost minimization (e.g., as in \cite{yang2014tight}, \cite{senjyu2003fast}, or \cite{dieu2006enhanced}).

Finally, in Table \ref{tab:118bus_solcomp}, we compare again the commitment solutions for the UC-ACOPF and the UC-DCOPF problems from the literature, along with those obtained using our direct approaches. All commitments were selected from those presented in the previous tables with a penalty weight of $5.\mathrm{e}{6}$, which also returned a feasible point for the UC-ACOPF problem. The rescale-and-round strategy over-committed plants compared to the solutions obtained using the optimization approaches over the binary decision variables in \cite{castillo2016unit} and \cite{liu2018global}. Decommitment heuristics, such as those in \cite{borghetti2003lagrangian}, could be applied, although maintaining AC-feasibility may prove difficult for highly constrained power grids. Our least costly solution was obtained using Re-Power rescaling and the NR formula; however, this method did not prove robust in terms of feasibility across the numerical experiments. Compared to the local solution in \cite{castillo2016unit}, our cost overruns were less than $10000$\$, while the underestimated cost of solving the (UC-DCOPF) problem was over $20000$\$. Although our direct approach is not a full optimization heuristic, it could be a good alternative for fast applications where obtaining an AC-feasible point is crucial.
\begin{table}
\centering
\footnotesize
\begin{tabular}{|c|c|c|c|cc|ccc|}
\hline
 & (UC-DCOPF) & \multicolumn{2}{|c|}{Literature (UC-ACOPF)} & \multicolumn{5}{|c|}{Direct strategies} \\
\hline
& \cite{castillo2016unit} & \cite{liu2018global} & \cite{castillo2016unit}  & \multicolumn{2}{|c|}{NR} &  \multicolumn{3}{|c|}{UC-ER} \\
& & Global & Local & Re-RUC & Re-Power & None & Re-RUC & Re-Power \\
\hline\hline
 Cost($\$$) & 814715 & 835926 & 843591  & 847740  &  846543 &  850817 & 851456 & 848837\\
 \hline
g.u. 1  & -    & -         & -     & -      & -      & -        & -    & -    \\
g.u. 2  & -    & -         & -     & -      & -      & -        & -    & -       \\
g.u. 3  & -    & -         & -     & 1-24   & 9-23   & 1-24     & 1-24 & 9-10,18  \\
g.u. 4  & 1-24 & 1-10,24   & 1-24  & 1-24   & 1-24   & 1-24     & 1-24 & 1-24    \\
g.u. 5  & 1-24 & 1-24      & 1-24  & 1-24   & 1-24   & 1-24     & 1-24 & 1-24    \\
g.u. 6  & -    & -         & -     & -      & -      & -        & -    & -       \\
g.u. 7  & -    & 11-22     & 10-23 & 6-24   & 6-24   & 7-23     & 7-23 & 7-23    \\
g.u. 8  & -    & -         & -     & -      & -      & -        & -    & -       \\
g.u. 9  & -    & -         & -     & -      & -      & -        & -    & -       \\
g.u. 10 & 1-24 & 1-2,12-24 & 1-24  & 1-24   & 1-24   & 1-24     & 1-24 & 1-24    \\
g.u. 11 & 1-24 & 1-24      & 1-24  & 1-24   & 1-24   & 1-24     & 1-24 & 1-24    \\
g.u. 12 & -    & -         & -     & -      & -      & -        & -    & -       \\
g.u. 13 & -    & -         & -     & -      & -      & -        & -    & -       \\
g.u. 14 & -    & 10-22     & 9-22  & 8-24   & 8-24   & 8-24     & 8-24 & 8-24    \\
g.u. 15 & -    & -         & -     & -      & -      & -        & -    & -       \\
g.u. 16 & 10-18& 9-16      & 10-23 & 8-24   & 8-24   & 8-24     & 8-23 & 8-23   \\
g.u. 17 & -    & -         & -     & -      & -      & -        & -    & -       \\
g.u. 18 & -    & -         & -     & -      & -      & -        & -    & -       \\
g.u. 19 & -    & -         & -     & 10-22  & 10-22  & 9-14     & 8-15 & 8-15    \\
g.u. 20 & 1-24 & 1-24      & 1-24  & 1-24   & 1-24   & 1-24     & 1-24 & 1-24    \\
g.u. 21 & -    & 8-24      & 8-24  & 1-24   & 1-24   & 1-24     & 1-24 & 1-24    \\
g.u. 22 & -    & -         & -     & 9-23   & 9-23   & 9-19     & 9-22 & 9-22    \\
g.u. 23 & -    & -         & -     & -      & -      & -        & -    & -       \\
g.u. 24 & -    & 9-23      & 9-23  & 7-24   & 7-24   & 7-24     & 7-24 & 7-24    \\
g.u. 25 & 10-22& -         & -     & 1-24   & 1-24   & 1-24     & 1-23 & 1-23    \\
g.u. 26 & -    & -         & -     & -      & -      & -        & -    & -       \\
g.u. 27 & 1-24 & 1-2,13-24 & 1-24  & 1-24   & 1-24   & 1-24     & 1-24 & 1-24    \\
g.u. 28 & 1-24 & 1-24      & 1-24  & 1-24   & 1-24   & 1-18     & 1-24 & 1-24    \\
g.u. 29 & 1-24 & 1-24      & 1-24  & 1-24   & 1-24   & 1-24     & 1-24 & 1-24    \\
g.u. 30 & 8-23 & 1-24      & 7-24  & 1-24   & 1-24   & 1-24     & 1-24 & 1-24    \\
g.u. 31 & -    & -         & -     & -      & -      & -        & -    & -       \\
g.u. 32 & -    & -         & -     & -      & -      & -        & -    & -       \\
g.u. 33 & -    & -         & -     & -      & -      & -        & -    & -       \\
g.u. 34 & 7-24 & 7-24      & 7-24  & 1-24   & 1-24   & 1-24     & 1-24 & 1-24    \\
g.u. 35 & 1-24 & 1-24      & 1-24  & 1-24   & 1-24   & 1-24     & 1-24 & 1-24    \\
g.u. 36 & 1-24 & 1-24      & 1-24  & 1-24   & 1-24   & 1-24     & 1-24 & 1-24    \\
g.u. 37 & 1-24 & 8-23      & 11-17 & 1-24   & 1-24   & 1-24     & 1-24 & 1-24    \\
g.u. 38 & -    & -         & -     &8,17-28 & 8,17-18& 8-9,18   &8-9,18& 8-10,18  \\
g.u. 39 & 1-24 & -         & 1-24  & -      & -      & 18-24    & 18-24& 18-24   \\
g.u. 40 & 1-24 &1-10,22-24 & 1-24  & 1-24   & 1-24   & 1-24     & 1-24 & 1-24    \\
g.u. 41 & -    & -         & -     & -      & -      & -        & -    & -       \\
g.u. 42 & -    & -         & -     & -      & -      & -        & -    & -       \\
g.u. 43 & 9-24 & 1-24      & 10-23 & 1-24   & 1-24   & 1-24     & 1-24 & 1-24    \\
g.u. 44 & -    & -         & -     & -      & -      & 10-17    & 10-17& 10-17   \\
g.u. 45 & 1-24 & 1-24      & 1-24  & 1-24   & 1-24   & 1-24     & 1-24 & 1-24    \\
g.u. 46 & -    & -         & -     & -      & -      & -        & -    & -       \\
g.u. 47 & -    & -         & -     & -      & -      & -        & -    & -       \\
g.u. 48 & -    & -         & -     & 8-23   & 8-23   & 8-23     & 8-23 & 8-23    \\
g.u. 49 & -    & -         & -     & -      & -      & -        & -    & -       \\
g.u. 50 & -    & -         & -     & -      & -      & -        & -    & -       \\
g.u. 51 & -    & 9-13      & 9-24  & 8-23   & 8-23   &8-12,18-22& 8-22 & 8-22    \\
g.u. 52 & -    & 14-23     & -     & 8-24   & 8-24   & 8-23     & 8-23 & 8-23    \\
g.u. 53 & 8-24 & 7-24      & 1-24  & 1-24   & 1-24   & 1-24     & 1-24 & 1-24    \\
g.u. 54 & 9-23 & 9-23      & 9-20  & 9-23   & 9-23   & 8-23     & 9-23 & 9-23    \\
 \hline
\end{tabular}
\caption{Comparison between commitment and cost solutions for the 118-bus test system example.}
\label{tab:118bus_solcomp}
\end{table}
\normalsize
%
\section{Integration of the relax-and-round strategies into a FP algorithm}\label{Sect5}
In this section, we evaluate the broad applicability of the previously presented relax-and-round strategies in combination with other heuristics from the literature by integrating them into a Feasibility Pump algorithm.

The Feasibility Pump (FP) method is a heuristic for finding feasible points for MILP problems, first proposed in \cite{fischetti2005feasibility}. The core idea of this method is to relax the integrality condition, solve the resulting relaxed problem, and then iteratively solve a family of LP problems until an integer feasible point is found. The formulation of the LP problems retains the original constraints but changes the objective function to be the $L_1$-norm of the difference between the relaxed solution and its integer-feasible rounded counterpart. To improve the quality of the solutions returned by this heuristic, the so-called Objective Feasibility Pump (OFP) was proposed in \cite{achterberg2007improving}. For an introduction to the method, we refer the reader to the papers cited above. Here, we have implemented an Objective Feasibility Pump algorithm for the Unit Commitment problem with ACPF constraints. More precisely, we extend the work in \cite{li2024objective} by considering the minimum up/down constraints in the UC formulation, and by broadening the time horizon to 24 hours. In our case, the FP method is implemented as a heuristic to find feasible points for a MINLP problem with a strong combinatorial nature, involving the solution of a family of non-convex NLP problems. We have implemented the ideas behind the OFP algorithm for the UC-ACOPF problem in a straightforward way. A rigorous adaptation of the FP heuristic to this analytical framework is beyond the scope of this paper; we refer the reader to \cite{d2012storm} for a discussion on extending the FP heuristic to general MINLP problems. The main goal is to use the rescaling strategy and the rounding formulas presented in the previous section together with our implementation of the FP method to check if the performance of the heuristic improves, compared to using the standard approach in the literature.

Let us consider the vector $\boldsymbol{u}^{\text{in}}$, consisting of binary values $u_{g,t}^{\text{in}}$ for each $g$ in the set of thermal generators $G_{\text{th}}$ and each time period $t\in T$. Each value $u_{g,t}^{\text{in}}$ corresponds specifically to the decision variable $\tu_{g,t}$ of the problem RUC-ACOPF, grouped in a vector denoted by $\boldsymbol{\tu}$. The $L1$-norm of the difference between $\boldsymbol{u}^{\text{in}}$ and $\boldsymbol{\tu}$ can be computed as
$$\| \boldsymbol{u}^{\text{in}}-\boldsymbol{\tu}\|_{L_1}= {\sum_{\left(g,t\right)\in G_{\text{th}}\times T} \left|u_{g,t}^{\text{in}}-\tu_{g,t}\right|}= \sum_{t\in T}\sum_{g\in G_{\text{th}}} \left|u_{g,t}^{\text{in}}-\tu_{g,t}\right|.$$
Since $\boldsymbol{u}^{\text{in}}$ is a vector of binary variables, the previous formula can also be expressed as

$$\| \boldsymbol{u}^{\text{in}}-\boldsymbol{\tu}\|_{L_1} = \sum_{\left(g,t\right)\in G_{\text{th}}\times T: u_{g,t}^{\text{in}}=0} \tu_{g,t} + \sum_{\left(g,t\right)\in G_{\text{th}}\times T: u_{g,t}^{\text{in}}=1} \left(1- \tu_{g,t}\right).$$
The cost function for the sequence of NLP problems iteratively solved in the objective feasibility pump algorithm is defined through a convex combination of the UC cost function $f_{\text{UC}}$ and the  $L_1$-norm difference between the integer and relaxed decision variables. Thus:
\begin{equation}
\label{FP_cost}
f_{\text{FP}}\left(\boldsymbol{x},\boldsymbol{y}\right) \defeq \frac{1-\alpha}{W_{L_1}}\|\boldsymbol{u}^{in}-\boldsymbol{\tu}\|_{L_1} + \frac{\alpha}{W_{\text{UC}}}f_{\text{UC}}\left(\boldsymbol{x},\boldsymbol{y}\right),
\end{equation}
where $W_{L_1}$ and $W_{\text{UC}}$ are weights used to normalize the difference in magnitude between the two terms of the cost function, and the convex combination factor $\alpha \in \left(0,1\right)$ is updated after each iteration.

Let us consider the set of  NLP problems denoted as (FP-RUC), parameterized by $\boldsymbol{u}^{in}$ and $\alpha$:
\begin{equation}
\left(\text{FP-RUC}\right)\left(\boldsymbol{u}^{in},\alpha\right) \quad
\left\{
\begin{array}{rl}
    \min & \left(\ref{FP_cost}\right), \\[0.1ex]
   \text{s.t.:} & \left(\ref{UC_BVL}\right)-\left(\ref{ACPF_PV_LineLim}\right), \\[0.1ex]
    & V_{i,t}\in\left[V_{\min,i},V_{\max,i}\right],\ \theta_{i,t}\in\left[-\pi,\pi\right],\ i\in I, t\in T,\\[0.1ex]
    & \theta_{i,t} = \theta_{i,0},\ \forall t\in T,\ \text{ for the reference node},\\[0.1ex]
    &   Q^G_{g,t}\in\left[Q_{\min,g},Q_{\max,g}\right],\ \forall \left(g,t\right)\in G_{\text{sc}}\times T,\\[0.1ex]
    & P^G_{g,t}\in\left[0,P_{\max,g}\right],\ \forall \left(g,t\right)\in G_{\text{th}}\times T,\\
    &   Q^G_{g,t}\in\left[\min\{Q_{\min,g},0\},\max\{Q_{\max,g},0\}\right],\ \forall \left(g,t\right)\in G_{\text{th}}\times T,\\[0.1ex]
    &P_{g,t}^{\text{res}}\in\left[0,P_{\max,g}- P_{\min,g}\right],\ \forall \left(g,t\right)\in G_{\text{th}}\times T, \\[0.1ex]
    & u_{g,t}, v_{g,t}, w_{g,t} \in \left[0,1\right],\ \forall \left(g,t\right)\in G_{\text{th}}\times T.\\[0.1ex]
\end{array}
\right.
\end{equation}
Each distinct pair of values for $\left(\boldsymbol{u}^{in},\alpha\right)$ defines a unique problem instance within this family. Our relax-and-round approach using the FP heuristic is described in Algorithms \ref{alg:relax_fp} and \ref{alg:fp}.

\begin{algorithm}[H]
\small
\caption{Relax-and-round with Feasibility Pump}\label{alg:relax_fp}
\begin{algorithmic}[1]
\State Consider (UC-ACOPF).
\State \text{Step 1}: Solve (RUC-ACOPF). Let $\tilde{u}_{g,t}$ be the relaxed commitment.
\State \text{Step 2}: Round $\tilde{u}_{g,t}$ values to obtain $\boldsymbol{u}^{in}$
\State{\text{Step 3}: Apply Algorithm \ref{alg:fp} implementing the Feasibility Pump using $\boldsymbol{u}^{in}$ as the initial integer vector.}
\end{algorithmic}
\end{algorithm}
\begin{algorithm}
\caption{Feasibility Pump algorithm}\label{alg:fp}
\begin{algorithmic}[1]
\Require (UC-ACOPF) problem data, $\boldsymbol{u}^{\text{in}}$ initial integer value, and parameters: \textit{maxit}, \textit{maxrst}, $\alpha_0$, $\varphi_{\alpha}$, $\delta_{\alpha}$ and $S$.
\State Initialize: $\textit{rst}=0$, $\boldsymbol{u}_0^{\text{in}}=\boldsymbol{u}^{\text{in}}$.
\For{$k \gets 1$ \textbf{to} $\textit{maxit}$}
    \State Solve (FP-RUC)$\left(\boldsymbol{u}_{k-1}^{\text{in}}, \alpha_{k-1}\right)$. Let $\left(\boldsymbol{\hat{x}}_k,\boldsymbol{\hat{y}}_k\right)$ be the solution, and $\boldsymbol{\hat{u}}_k$ the relaxed commitment.
    \If{ $\boldsymbol{\hat{u}}_k$ is integer}
        \State \textbf{Return} $\left(\boldsymbol{\hat{x}}_k,\boldsymbol{\hat{y}}_k\right)$ and \textbf{End}.
    \Else
        \State Round $\boldsymbol{\hat{u}}_k$ to obtain $\boldsymbol{u}_k^{\text{in}}$.
    \EndIf
    \If{$\left(\boldsymbol{u}_k^{\text{in}},\alpha_k\right)$ is a stationary point}
        \State Apply a flip to the binary values of $S$ random components of $\boldsymbol{u}_k^{\text{in}}$, and $\textit{rst}=\textit{rst}+1$.
        \State \textbf{if} $\textit{rst}>\textit{maxrst}$ \textbf{Exit for, go to} \textbf{15:}
    \EndIf
    \State Update: $\alpha_{k}=\alpha_{k-1}\varphi_{\alpha}$.
\EndFor
\State Choose $\left(\boldsymbol{\hat{x}}_h,\boldsymbol{\hat{y}}_h\right)$, $h\in\{1,\dots,k\}$, with closest  $\boldsymbol{\hat{u}}_{h}$ to its projection in \textbf{7:} in the $L_1$-norm sense.
\State Round $\boldsymbol{\hat{u}}_{h}$ and construct $\boldsymbol{{y}}_h^{\text{in}}$.
\State Fix $\boldsymbol{{y}}_h^{\text{in}}$ in (UC-ACOPF), and solve the resulting problem (ACOPF). Let $\boldsymbol{x}_{\text{AC}}$ be the solution.
\State \textbf{Return} $\left(\boldsymbol{x}_{\text{AC}},\boldsymbol{{y}}_h^{\textit{in}}\right)$ and \textbf{End}.
\end{algorithmic}
\end{algorithm}

Let us now discuss the specifics of Algorithm \ref{alg:fp}. After solving (FP-RUC)$\left(\boldsymbol{u}_{k-1}^{\text{in}},\alpha_{k-1}\right)$ at iteration $k$, we denote by $\boldsymbol{\hat{u}}_k$ the values for the relaxed unit commitment in the solution. In lines 4-8 of the algorithm, if $\boldsymbol{\hat{u}}_k$ is an integer vector, we declare that the pair $\left(\boldsymbol{\hat{x}}_k,\boldsymbol{\hat{y}}_k\right)$ is a feasible point for the UC-ACOPF problem and the algorithm finishes; otherwise, we round $\boldsymbol{\hat{u}}_k$ to obtain the next candidate for the integer commitment, $\boldsymbol{u}_k^{\text{in}}$. Notice that the solutions of problems (FP-RUC)$\left(\boldsymbol{u}_{k-1}^{\text{in}},\alpha_{k-1}\right)$ can converge to the integer vector $\boldsymbol{u}_{k-1}^{\text{in}}$ only if it satisfies the combinatorial constraints of the UC problem. If the newly obtained integer vector $\boldsymbol{u}_k^{\text{in}}$ matches any previously obtained integer solution, say $\boldsymbol{u}_{\tilde{k}}^{\text{in}}$ for some earlier iteration $\tilde{k}<k$ while also satisfying $\alpha_{\tilde{k}}-\alpha_{k}\le \delta_{\alpha}$ we consider the pair $\left(\boldsymbol{u}_k^{\text{in}},\alpha_k\right)$ to be a stationary point for the algorithm. To avoid cyclic behavior, after rounding $\boldsymbol{\hat{u}}_k$, we verify if a stationary point is reached, in which case we apply a perturbation to $\boldsymbol{u}_k^{\text{in}}$ by flipping the binary values of $S$ random components. This perturbation mechanism can be seen as a `reset', as we also update $\alpha_k=\alpha_0/\left(\textit{rst}+1\right)$, where \textit{rst} counts the number of resets, thereby partially restarting the algorithm. To prevent excessive resets, we set a maximum number, \textit{maxrst}, as an early stopping criterion. The rounding procedure described in Algorithm \ref{alg:relax_fp} (Step 2) and Algorithm \ref{alg:fp} (lines 7 and 16) are independent. This means any rescaling and rounding formula combinations, as presented in the previous sections, can be applied. The naive rounding formula with no rescale (None$+$NR) projects the relaxed solution to the nearest integer vector in the $L_1$-norm sense. However, not all these integer combinations are feasible under the UC problem constraints, highlighting the importance of combining these techniques effectively for better solution feasibility and quality. In this sense, we expect that our proposed rescaling$+$rounding formulas provide a more accurate direct approximation of the projection of $\boldsymbol{\hat{u}}_k$ in the UC skeleton solution space.

In the following, we present the results of the FP algorithm applied to the test system introduced in the previous sections, i.e., the 6-bus and the 118-bus test systems. To solve the non-convex NLP problems, we once again employed the PSLP method. For both problems (RUC+ACOPF) and (ACOPF), a penalty weight of $5.\mathrm{e}{6}$ was utilized, whereas for problem (FP-RUC), the penalty weight was set to $500$. Furthermore, the parameters for the cost function were defined as follows: we set $W_{L_1} = t_f\left|G_{\text{th}}\right|$ and $W_{\text{UC}} = t_f\left|G_{\text{th}}\right|\max_{g\in G_{\text{th}}}\{a_{g,2}\}$, where $\left|G_{\text{th}}\right|$ denotes the number of thermal generators in the system.

For the FP algorithm, we used the following parameter settings: $\textit{maxrst}=3$, $\alpha_0=0.75$, $\varphi_{\alpha}=0.5$, and $\varphi_{\beta}=0.075$. The parameter $S$ was chosen as a percentage of $t_f\left|G_{\text{th}}\right|$, and the results for values corresponding to its $10\%$, $25\%$, $50\%$, $75\%$, and $100\%$ were compared. The value $100\%$ corresponds to the total number of relaxed variables that must be rounded for the problem. The larger the value of $S$, the more we rely on random perturbations as a multi-start mechanism. Notice that the random flips introduce stochasticity into the results. As we will show, different values for $S$ yield varying performance depending on the combination of rescaling approach and rounding formula used.

The results summarizing the performance of the FP algorithm for both test systems are shown in Tables \ref{tab:6bus_summary_FP} and \ref{tab:118bus_summary_FP}. They are based on 10 runs per configuration, with the statistics representing the range of values obtained in these experiments. In both test cases, Step 2 of Algorithm \ref{alg:relax_fp} employs the NR formula without rescaling, thus initializing the FP algorithm with a non-feasible integer solution (recall Tables \ref{tab:6bus_summary} and \ref{tab:118bus_summary_naive}).

\begin{table}[!t]
\centering
\begin{tabular}{lcccccccc}
\toprule
\multicolumn{9}{c}{ \textbf{6-bus test system}} \\
\toprule
\multicolumn{4}{c}{Feasibility Pump} & \multicolumn{1}{c}{PSLPs} &  \multicolumn{4}{c}{Results} \\
    \midrule
  \text{Rescale$+$Round} & $S$  &  Iter. & $rst$ & Total iter. &  Feas. & Cost($\$$) & Time (s) & $\%$ \\
\toprule
None+NR & $7$ &  13-15 & 3  &  $312-377$ & \xmark & $94704-99641$ & $38-48$ & $100\%$\\
 & $18$ &  15 & 3  &  $397-412$ & \xmark & $94704-102115$ & $48-51$ & $100\%$\\
 & $36$ &  15 & 3  &  $394-423$ & \xmark & $94704$ & $47-51$ & $90\%$\\
 & $36$ & 12 & 2  &  $288$ & \cmark & $105701$ & $35$ & $10\%$ \\
 & $54$ &  15-16 & 3  &  $398-447$ & \xmark & $94704-103629$ & $50-54$ & $90\%$\\
 & $54$ & 11 & 2  &  $266$ & \cmark & $107801$ & $32$ & $10\%$ \\
 & $72$ &  15 & 3  &  $394-418$ & \xmark & $94704-104999$ & $47-51$ & $100\%$\\
\hline
Re-RUC+NR & $7$ &  15 & 3  & $427-476$ & \xmark & $101796-103612$ & $55-56$ & $20\%$\\
          & $7$ &  7-16 & 1-3  & $186-470$ & \cmark & $102485-103516$ & $22-57$ & $80\%$\\
 & $18$ &  15 & 3  & $465$ & \xmark & $101058$ & $55$ & $10\%$\\
 & $18$ &  7-11 & 1-2  & $184-313$ & \cmark & $103264-105640$ & $22-38$ & $90\%$\\
 & $36$ &  7 & 1  & $186-190$ & \cmark & $104014-106166$ & $22-23$ & $100\%$\\
 & $54$ &  7-11 & 1-2  & $183-319$ & \cmark & $106240-107439$ & $22-38$ & $100\%$\\
 & $72$ &  7 & 1  & $184-191$ & \cmark & $106698-108507$ & $22-23$ & $100\%$\\
\hline

Re-RUC+UC-ER & Any & 2 & 0  & $37$ & \cmark & $102146$ & $5$ & $100\%$\\

\bottomrule
\end{tabular}
    \caption{Summary of the FP algorithm performance for the 6-bus test system example.}
    \label{tab:6bus_summary_FP}
\end{table}

For the 6-bus test system, the None$+$NR rounding approach produced infeasible solutions in $90\%$ of the cases. Notably, for the remaining $10\%$, the algorithm required two restarts to obtain a feasible solution. In contrast, when using Re-RUC$+$NR with $S=t_f\left|G_{\text{th}}\right|$, a feasible solution was found after only one restart, yielding a different feasible solution in each experiment. As the value of $S$ decreased, the Re-RUC$+$NR approach had increasing difficulty in identifying feasible solutions. Finally, the Re-RUC$+$UC-ER method achieved convergence in only two iterations, eliminating the need for any random perturbations.

For the 118-bus test system, we report only the feasible points obtained for the different values of $S$ and for each combination of rescaling and rounding strategies. The increased size of this test system negatively affected the numerical performance of the FP algorithm, making it challenging to fine-tune the parameters for both the PSLP and FP algorithms. This issue should be addressed in future work. For completeness, we adopted a higher feasibility threshold of $1.\mathrm{e}{-5}$ for this test system and verified solution feasibility using our Fortran code. The None$+$NR approach did not consistently yield feasible solutions for any value of $S$, except for a single feasible solution found for $S=t_f\left|G_{\text{th}}\right|0.5$. The Re-RUC$+$NR and Re-RUC$+$UC-ER approaches showed contrasting performance trends based on the value of $S$. For Re-RUC$+$UC-ER, using $S=t_f\left|G_{\text{th}}\right|$ yielded the best results, achieving a feasibility rate of $100\%$ within the set tolerance. As $S$ decreased, the feasibility rate for Re-RUC$+$UC-ER gradually decreased to $30\%$ at $S=t_f\left|G_{\text{th}}\right|0.25$. Conversely, Re-RUC$+$NR performed better for smaller values of $S$, also achieving $30\%$ feasibility at $S=t_f\left|G_{\text{th}}\right|0.25$, but showing reduced performance for other $S$. For smaller values of $S$, none of the rescaling and rounding combinations produced feasible points. This behavior aligns with previous findings on the usefulness of a multi-start mechanism in FP implementations for MINLP problems, as reported in \cite{d2012storm}.

\begin{table}
\centering
\begin{tabular}{lccccccc}
\toprule
\multicolumn{7}{c}{ \textbf{118-bus test system}} \\
\toprule
\multicolumn{2}{c}{Feasibility Pump} & \multicolumn{1}{c}{PSLPs} & \multicolumn{2}{c}{Feasibility} &  \multicolumn{3}{c}{Results} \\
    \midrule
  \text{Rescale$+$Round} & $S$  & Total iter. & AC & UC & Cost($\$$) & Time (s) & $\%$ \\
\toprule
None+NR & $324$   & $799$ & $1.\mathrm{e}{-5}$ & \cmark & $929131$ & $4188$ & $10\%$\\
\hline
Re-RUC+NR & $130$   & $1006$ & $1.\mathrm{e}{-5}$ & \cmark & $871041-871746$ & $3524$ & $10\%$\\

 & $324$   & $520-560$ & $1.\mathrm{e}{-5}$ & \cmark & $871041-871746$ & $2204-2320$ & $30\%$\\

 & $648$   & $3817$ & $1.\mathrm{e}{-5}$ & \cmark & $907184$ & $15381$ & $10\%$\\
\hline
 Re-RUC+UC-ER & $324$   & $720-2206$ & $1.\mathrm{e}{-5}$ & \cmark & $873450-892570$ & $3177-6554$ & $30\%$\\
 & $648$   & $2125-6367$ & $1.\mathrm{e}{-5}$ & \cmark & $889443-911307$ & $5721-14747$ & $40\%$\\
 & $972$   & $400-2084$ & $1.\mathrm{e}{-5}$ & \cmark & $899136-920748$ & $1754-14391$ & $60\%$\\
 & $1296$   & $360-2620$ & $1.\mathrm{e}{-5}$ & \cmark & $899924-929276$ & $1584-3120$ & $100\%$\\

\bottomrule
\end{tabular}
    \caption{Summary of the FP algorithm performance for the 118-bus test system example.}
    \label{tab:118bus_summary_FP}
\end{table}
\section{Conclusions}
Finding solutions, or even feasible points, for the UC problem with ACPF constraints remains a challenging task. Many heuristics in the literature rely on relaxing the integrality of the commitment variables, followed by rounding to restore integer feasibility. In this work, we proposed two novel heuristics designed to enhance relax-and-round strategies. The first heuristic applies a physics-based rescaling to the relaxed commitment variables before rounding, while the second offers an improved rounding formula that guarantees combinatorial feasibility and preserves active power linear balance in the system.

Both heuristics were evaluated using a direct relax-and-round approach on the 6-bus and 118-bus test systems. This approach shows promise as an effective alternative for quickly obtaining AC-feasible solutions in some cases. The incorporation of these heuristics into the Feasibility Pump algorithm demonstrated broader applicability in combination with other strategies from the literature. Compared to the standard approach —for which no rescaling is applied and the rounding is performed around the $0.5$ threshold— our proposed heuristics consistently yielded more feasible points for the UC-ACOPF problem and improved the convergence rates of the strategies tested in this work.

The non-convex NLPs involved in both algorithms were solved locally using the PSLP method. The reported solving times are mainly driven by the convergence behavior of the PSLP algorithm. Future work could enhance the performance of these methods by employing more robust or ad hoc solvers for the ACOPF problems. Additionally, exploring global optimization techniques could be a promising direction for further research.

\section*{Acknowledgments}
 This work was partially supported by MCIN /AEI /10.13039/501100011033 / FEDER, UE through grant PID2021-122625OB-I00,  by  Ministerio de Ciencia, Innovación y Universidades through the Plan Nacional de I+D+i (MTM2017-86459-R) and the grant PRE2018-083893, and by Xunta de Galicia funds through grant GRC GI-1563 - ED431C 2021/15 and. It has also been supported by the German Research Foundation (DFG) under grant GO 1920/11-1 and 12-1.

\bibliography{references}

\begin{thebibliography}{10}

\bibitem{garver1962power}
Len~L Garver.
\newblock Power generation scheduling by integer programming-development of theory.
\newblock {\em Transactions of the American Institute of Electrical Engineers. Part III: Power Apparatus and Systems}, 81(3):730--734, 1962.

\bibitem{padhy2004unit}
Narayana~Prasad Padhy.
\newblock Unit commitment-a bibliographical survey.
\newblock {\em IEEE Transactions on Power Systems}, 19(2):1196--1205, 2004.

\bibitem{li2024objective}
Peijie Li, Jianming Su, and Xiaoqing Bai.
\newblock An objective feasibility pump method for optimal power flow with unit commitment variables.
\newblock {\em Electric Power Systems Research}, 236:110928, 2024.

\bibitem{wang2013stochastic}
Jiadong Wang, Jianhui Wang, Cong Liu, and Juan~P Ruiz.
\newblock Stochastic unit commitment with sub-hourly dispatch constraints.
\newblock {\em Applied Energy}, 105:418--422, 2013.

\bibitem{zhang2023mathematical}
Zihan Zhang, Mingbo Liu, Min Xie, and Ping Dong.
\newblock A mathematical programming--based heuristic for coordinated hydrothermal generator maintenance scheduling and long-term unit commitment.
\newblock {\em International Journal of Electrical Power \& Energy Systems}, 147:108833, 2023.

\bibitem{frangioni2008tighter}
Antonio Frangioni, Claudio Gentile, and Fabrizio Lacalandra.
\newblock Tighter approximated milp formulations for unit commitment problems.
\newblock {\em IEEE Transactions on Power Systems}, 24(1):105--113, 2009.

\bibitem{caroe1998two}
Claus~C. Car{\o}e and R{\"u}diger Schultz.
\newblock A two-stage stochastic program for unit commitment under uncertainty in a hydro-thermal power system.
\newblock {\em Preprint SC 98-11}, 1998.

\bibitem{wang2008security}
Jianhui Wang, Mohammad Shahidehpour, and Zuyi Li.
\newblock Security-constrained unit commitment with volatile wind power generation.
\newblock {\em IEEE Transactions on Power Systems}, 23(3):1319--1327, 2008.

\bibitem{wang2011wind}
J~Wang, A~Botterud, R~Bessa, H~Keko, L~Carvalho, D~Issicaba, J~Sumaili, and V~Miranda.
\newblock Wind power forecasting uncertainty and unit commitment.
\newblock {\em Applied Energy}, 88(11):4014--4023, 2011.

\bibitem{carrion2006computationally}
Miguel Carri{\'o}n and Jos{\'e}~M Arroyo.
\newblock A computationally efficient mixed-integer linear formulation for the thermal unit commitment problem.
\newblock {\em IEEE Transactions on Power Systems}, 21(3):1371--1378, 2006.

\bibitem{morales2013tight}
Germ{\'a}n Morales-Espa{\~n}a, Jesus~M Latorre, and Andres Ramos.
\newblock Tight and compact milp formulation for the thermal unit commitment problem.
\newblock {\em IEEE Transactions on Power Systems}, 28(4):4897--4908, 2013.

\bibitem{yang2021two}
Linfeng Yang, Wei Li, Yan Xu, Cuo Zhang, and Shifei Chen.
\newblock Two novel locally ideal three-period unit commitment formulations in power systems.
\newblock {\em Applied Energy}, 284:116081, 2021.

\bibitem{atakan2017state}
Semih Atakan, Guglielmo Lulli, and Suvrajeet Sen.
\newblock A state transition mip formulation for the unit commitment problem.
\newblock {\em IEEE Transactions on Power Systems}, 33(1):736--748, 2017.

\bibitem{chen2019modeling}
Honglin Chen, Mingbo Liu, Ying Cheng, and Shunjiang Lin.
\newblock Modeling of unit commitment with {AC} power flow constraints through semi-continuous variables.
\newblock {\em IEEE Access}, 7:52015--52023, 2019.

\bibitem{wu2010accelerating}
Lei Wu and Mohammad Shahidehpour.
\newblock Accelerating the benders decomposition for network-constrained unit commitment problems.
\newblock {\em Energy Systems}, 1(3):339--376, 2010.

\bibitem{castillo2016unit}
Anya Castillo, Carl Laird, C{\'e}sar~A Silva-Monroy, Jean-Paul Watson, and Richard~P O’Neill.
\newblock The unit commitment problem with {AC} optimal power flow constraints.
\newblock {\em IEEE Transactions on Power Systems}, 31(6):4853--4866, 2016.

\bibitem{holzer2024grid}
Jesse~T Holzer, Carleton~J Coffrin, Christopher DeMarco, Ray Duthu, Stephen~T Elbert, Brent~C Eldridge, Tarek Elgindy, Manuel Garcia, Scott~L Greene, Nongchao Guo, et~al.
\newblock Grid optimization competition challenge 3 problem formulation.
\newblock Technical report, Pacific Northwest National Laboratory (PNNL), Richland, WA (United States), 2024.

\bibitem{yang2014tight}
Linfeng Yang, Jinbao Jian, Yunan Zhu, and Zhaoyang Dong.
\newblock Tight relaxation method for unit commitment problem using reformulation and lift-and-project.
\newblock {\em IEEE Transactions on Power Systems}, 30(1):13--23, 2014.

\bibitem{ostrowski2011tight}
James Ostrowski, Miguel~F Anjos, and Anthony Vannelli.
\newblock Tight mixed integer linear programming formulations for the unit commitment problem.
\newblock {\em IEEE Transactions on Power Systems}, 27(1):39--46, 2011.

\bibitem{bai2009semi}
X~Bai and H~Wei.
\newblock Semi-definite programming-based method for security-constrained unit commitment with operational and optimal power flow constraints.
\newblock {\em IET Generation, Transmission \& Distribution}, 3(2):182--197, 2009.

\bibitem{murillo2000parallel}
Carlos~E Murillo-Sanchez and Robert~J Thomas.
\newblock Parallel processing implementation of the unit commitment problem with full {AC} power flow constraints.
\newblock In {\em Proceedings of the 33rd Annual Hawaii International Conference on System Sciences}, pages 9--pp. IEEE, 2000.

\bibitem{niknam2009new}
Taher Niknam, Amin Khodaei, and Farhad Fallahi.
\newblock A new decomposition approach for the thermal unit commitment problem.
\newblock {\em Applied Energy}, 86(9):1667--1674, 2009.

\bibitem{nick2015security}
Mostafa Nick, Omid Alizadeh-Mousavi, Rachid Cherkaoui, and Mario Paolone.
\newblock Security constrained unit commitment with dynamic thermal line rating.
\newblock {\em IEEE Transactions on Power Systems}, 31(3):2014--2025, 2015.

\bibitem{cobos2018network}
Noemi~G Cobos, Jos{\'e}~M Arroyo, Natalia Alguacil, and Alexandre Street.
\newblock Network-constrained unit commitment under significant wind penetration: A multistage robust approach with non-fixed recourse.
\newblock {\em Applied energy}, 232:489--503, 2018.

\bibitem{senjyu2003fast}
Tomonobu Senjyu, Kai Shimabukuro, Katsumi Uezato, and Toshihisa Funabashi.
\newblock A fast technique for unit commitment problem by extended priority list.
\newblock {\em IEEE Transactions on Power Systems}, 18(2):882--888, 2003.

\bibitem{xavier2021learning}
{\'A}linson~S Xavier, Feng Qiu, and Shabbir Ahmed.
\newblock Learning to solve large-scale security-constrained unit commitment problems.
\newblock {\em INFORMS Journal on Computing}, 33(2):739--756, 2021.

\bibitem{tejada2019unit2}
Diego~A Tejada-Arango, Sara Lumbreras, Pedro S{\'a}nchez-Mart{\'\i}n, and Andres Ramos.
\newblock Which unit-commitment formulation is best? {A comparison framework}.
\newblock {\em IEEE Transactions on Power Systems}, 35(4):2926--2936, 2019.

\bibitem{sauer2014uplift}
William Sauer.
\newblock Uplift in rto and iso markets.
\newblock {\em Federal Energy Regulatory Commission, Tech. Rep}, 2014.

\bibitem{bai2015decomposition}
Yang Bai, Haiwang Zhong, Qing Xia, Chongqing Kang, and Le~Xie.
\newblock A decomposition method for network-constrained unit commitment with {AC} power flow constraints.
\newblock {\em Energy}, 88:595--603, 2015.

\bibitem{liu2018global}
Jianfeng Liu, Carl~D Laird, Joseph~K Scott, Jean-Paul Watson, and Anya Castillo.
\newblock Global solution strategies for the network-constrained unit commitment problem with {AC} transmission constraints.
\newblock {\em IEEE Transactions on Power Systems}, 34(2):1139--1150, 2018.

\bibitem{tuncer2022misocp}
Deniz Tuncer and Burak Kocuk.
\newblock An misocp-based decomposition approach for the unit commitment problem with {AC} power.
\newblock {\em IEEE Transactions on Power Systems}, 38(4):3388--3400, 2022.

\bibitem{montero2022review}
Luis Montero, Antonio Bello, and Javier Reneses.
\newblock A review on the unit commitment problem: Approaches, techniques, and resolution methods.
\newblock {\em Energies}, 15(4):1296, 2022.

\bibitem{aharwar2023unit}
Ankit Aharwar, Ram Naresh, Veena Sharma, and Vineet Kumar.
\newblock Unit commitment problem for transmission system, models and approaches: A review.
\newblock {\em Electric Power Systems Research}, 223:109671, 2023.

\bibitem{wuijts2024effect}
Rogier~Hans Wuijts, Marjan van~den Akker, and Machteld van~den Broek.
\newblock Effect of modelling choices in the unit commitment problem.
\newblock {\em Energy Systems}, 15(1):1--63, 2024.

\bibitem{fischetti2005feasibility}
Matteo Fischetti, Fred Glover, and Andrea Lodi.
\newblock The feasibility pump.
\newblock {\em Mathematical Programming}, 104:91--104, 2005.

\bibitem{lazic2016variable}
Jasmina Lazi{\'c}.
\newblock Variable and single neighbourhood diving for mip feasibility.
\newblock {\em Yugoslav Journal of Operations Research}, 26(2), 2016.

\bibitem{ma2020unit}
Ziming Ma, Haiwang Zhong, Qing Xia, Chongqing Kang, Qiang Wang, and Xin Cao.
\newblock A unit commitment algorithm with relaxation-based neighborhood search and improved relaxation inducement.
\newblock {\em IEEE Transactions on Power Systems}, 35(5):3800--3809, 2020.

\bibitem{gao4718365topology}
Liqian Gao, Jiakun Fang, Xiaomeng Ai, Lishen Wei, Shichang Cui, Wei Yao, and Jinyu Wen.
\newblock A topology-guided learning framework for security-constraint unit commitment.
\newblock {\em Available at SSRN 4718365}, 2024.

\bibitem{ordoudis2015stochastic}
Christos Ordoudis, Pierre Pinson, Marco Zugno, and Juan~M Morales.
\newblock Stochastic unit commitment via progressive hedging—extensive analysis of solution methods.
\newblock In {\em 2015 IEEE Eindhoven PowerTech}, pages 1--6. IEEE, 2015.

\bibitem{wu2023synergistic}
Jianghua Wu, Peter~B Luh, Yonghong Chen, Bing Yan, and Mikhail~A Bragin.
\newblock Synergistic integration of machine learning and mathematical optimization for unit commitment.
\newblock {\em IEEE Transactions on Power Systems}, 39(1):391--401, 2023.

\bibitem{kjeldsen2012heuristic}
Niels~Hvidberg Kjeldsen and Marco Chiarandini.
\newblock Heuristic solutions to the long-term unit commitment problem with cogeneration plants.
\newblock {\em Computers \& Operations Research}, 39(2):269--282, 2012.

\bibitem{zimmerman2024matpower}
Ray~D Zimmerman and C~Murillo-S{\'a}nchez.
\newblock Matpower user’s manual. 2024.
\newblock {\em URL https://matpower. org/docs/MATPOWER-manual. pdf. Accessed}, 27, 2020.

\bibitem{constante2022ac}
Gonzalo~E Constante-Flores, Antonio~J Conejo, and Feng Qiu.
\newblock {AC} network-constrained unit commitment via conic relaxation and convex programming.
\newblock {\em International Journal of Electrical Power \& Energy Systems}, 134:107364, 2022.

\bibitem{he2016robust}
Chuan He, Lei Wu, Tianqi Liu, and Mohammad Shahidehpour.
\newblock Robust co-optimization scheduling of electricity and natural gas systems via admm.
\newblock {\em IEEE Transactions on Sustainable Energy}, 8(2):658--670, 2016.

\bibitem{wang2021optimal}
Zhuo Wang, Luhao Wang, Zhengmao Li, Xingong Cheng, and Qiqiang Li.
\newblock Optimal distributed transaction of multiple microgrids in grid-connected and islanded modes considering unit commitment scheme.
\newblock {\em International Journal of Electrical Power \& Energy Systems}, 132:107146, 2021.

\bibitem{fu2006ac}
Yong Fu, Mohammad Shahidehpour, and Zuyi Li.
\newblock {AC} contingency dispatch based on security-constrained unit commitment.
\newblock {\em IEEE Transactions on Power Systems}, 21(2):897--908, 2006.

\bibitem{sawa2007security}
Toshiyuki Sawa, Yasuo Sato, Mitsuo Tsurugai, and Tsukasa Onishi.
\newblock Security constrained integrated unit commitment using quadratic programming.
\newblock In {\em 2007 IEEE Lausanne Power Tech}, pages 1858--1863. IEEE, 2007.

\bibitem{bazaraa2006nonlinear}
Mokhtar~S Bazaraa, Hanif~D Sherali, and Chitharanjan~M Shetty.
\newblock {\em Nonlinear programming: theory and algorithms}.
\newblock John Wiley \& Sons, 2006.

\bibitem{gurobi}
{Gurobi Optimization, LLC}.
\newblock {Gurobi Optimizer Reference Manual}, 2023.

\bibitem{Achterberg2009}
Tobias Achterberg.
\newblock {SCIP: solving constraint integer programs}.
\newblock {\em Mathematical Programming Computation}, 1(1):1--41, Jul 2009.

\bibitem{czyzyk_et_al_1998}
Joseph {Czyzyk}, Michael~P. {Mesnier}, and Jorge~J. {Mor{\'{e}}}.
\newblock The {NEOS} server.
\newblock {\em IEEE Journal on Computational Science and Engineering}, 5(3):68—75, 1998.

\bibitem{najman2019mccormick}
Jaromil Najman, Dominik Bongartz, Susanne Sass, Hatim Djelassi, Daniel Jungen, Wolfgang~R Huster, Jannik Burre, Kaan Karacasulu, Artur~M Schweidtmann, and Alexander Mitsos.
\newblock Mccormick-based algorithm for mixed-integer nonlinear global optimization.
\newblock In {\em 2019 AIChE Annual Meeting}. American Institute of Chemical Engineers, 2019.

\bibitem{fu2005security}
Yong Fu, Mohammad Shahidehpour, and Zuyi Li.
\newblock Security-constrained unit commitment with {AC} constraints.
\newblock {\em IEEE transactions on {P}ower {S}ystems}, 20(3):1538--1550, 2005.

\bibitem{tejada2019unit}
Diego~A Tejada-Arango, Sonja Wogrin, Pedro S{\'a}nchez-Mart{\i}n, and Andres Ramos.
\newblock Unit commitment with {ACOPF} constraints: Practical experience with solution techniques.
\newblock {\em 2019 IEEE Milan PowerTech}, pages 1--6, 2019.

\bibitem{dieu2006enhanced}
VN~Dieu and W~Ongsakul.
\newblock Enhanced augmented lagrangian hopfield network for unit commitment.
\newblock {\em IEE Proceedings-Generation, Transmission and Distribution}, 153(6):624--632, 2006.

\bibitem{borghetti2003lagrangian}
Alberto Borghetti, Antonio Frangioni, Fabrizio Lacalandra, and Carlo~Alberto Nucci.
\newblock Lagrangian heuristics based on disaggregated bundle methods for hydrothermal unit commitment.
\newblock {\em IEEE Transactions on Power Systems}, 18(1):313--323, 2003.

\bibitem{achterberg2007improving}
Tobias Achterberg and Timo Berthold.
\newblock Improving the feasibility pump.
\newblock {\em Discrete Optimization}, 4(1):77--86, 2007.

\bibitem{d2012storm}
Claudia D’Ambrosio, Antonio Frangioni, Leo Liberti, and Andrea Lodi.
\newblock A storm of feasibility pumps for nonconvex minlp.
\newblock {\em Mathematical Programming}, 136:375--402, 2012.

\end{thebibliography}

\end{document}